\pdfoutput=1
\documentclass{article}
\usepackage[utf8]{inputenc}
\usepackage{mathpazo} 
\usepackage{lmodern} 
\usepackage{subfig}
\usepackage{mathtools,amssymb}
\usepackage{graphicx,color}
\usepackage[version=3]{mhchem}
\usepackage{textcomp} 
\usepackage{amsmath,amsfonts,amssymb,epsfig,latexsym}
\usepackage{url}
\usepackage{graphics}
\usepackage{booktabs}
\usepackage{fancyhdr}
\usepackage[utf8]{inputenc}
\usepackage{float}
\usepackage{xcolor}
\usepackage{dsfont}
\newtheorem{tm}{Theorem}[section]
\newtheorem{df}[tm]{Definition}

\newtheorem{nap}[tm]{Remark}

\newcommand{\Z}{ \mathbb{Z} }

\newcommand{\R}{ \mathbb{R} }

\newcommand{\dist}{dist}
\usepackage[T1]{fontenc}
\usepackage[utf8]{inputenc}
\usepackage{bbm}
\usepackage{mathtools}
\DeclareMathOperator{\Var}{Var}

\title{Central Limit Theorem for Functionals of Persistence Diagrams in Germ-Grain Random Set Models with Applications to Goodness-of-Fit Testing}

\author{
Vesna Gotovac \DJ oga\v{s}\\
\small Faculty of Science, University of Split\\
\small Ru\dj era Bo\v{s}kovi\'ca 33, 21000 Split, Croatia\\
\small \texttt{vgotovac@pmfst.hr}
\and
Marcela Mandari\'c\\
\small Faculty of Science, University of Split\\
\small Ru\dj era Bo\v{s}kovi\'ca 33, 21000 Split, Croatia\\
\small \texttt{mmandaric@pmfst.hr}
}

\date{}

\begin{document}

\maketitle

\begin{abstract}
This paper establishes a central limit theorem (CLT) for functionals of
$M$-bounded persistence diagrams arising from germ--grain random set models.
Building on stabilisation methods for marked point processes, we show that,
under certain conditions, these topological summaries exhibit asymptotic
normality as the observation window increases, particularly for models with
exponential decay of correlations. These results are applied in
goodness-of-fit tests designed to detect spatial interactions such as
clustering or repulsion. Using test statistics derived from rectangular
partitions of persistence diagrams and functional summaries, such as the APF
or the support function of the lift zonoid, the study distinguishes between
different models. Finally, the methodology is applied to histological images
of breast tissue.
\end{abstract}

\noindent\textbf{Keywords:}
germ-grain model; persistence diagram; central limit theorem; exponential
decay of correlations

\medskip
\noindent\textbf{2020 Mathematics Subject Classification:}
60D05; 62R40; 55N31

\section{Introduction}

Random sets have gained significant attention in recent years as powerful tools for the statistical analysis of geometric structures arising in diverse scientific fields. Their versatility enables researchers to model complex natural phenomena such as spatial plant distributions in ecology \cite{moeller:2010}, the microscopic architecture of biological tissues in medicine \cite{hermann:2015}, and the structure of materials in physics and engineering \cite{neumann:2016, micheletti:2016}. A well-established theoretical framework supports the study of random sets \cite{matheron:1975, molchanov:2005, serra:1982}, making them an appealing foundation for developing new statistical methodologies.

In \cite{vesnamarcela}, a methodology for detecting outliers and testing the goodness-of-fit of random sets was introduced using tools from topological data analysis (TDA). This approach is very useful for investigating characteristics of random sets, such as their shapes or dependencies between components. These characteristics could indicate whether clustering or repulsion occurs. They can be revealed using topological data analysis, which can be used to investigate the topological properties of realisations of random sets. The core of any TDA method is the construction of a filtration and the tracking of the appearance or disappearance of topological features throughout the filtration. Persistent homology is used to track the persistence of topological features in the filtration as it grows. Homology groups were originally introduced to compare shapes by analysing their holes in different dimensions. The $0$th homology group, $H_0$, corresponds to the connected components of an object, while the $1$st homology group, $H_1$, corresponds to loops or one-dimensional holes. Usually, features that persist longer are more significant. Persistence diagrams (PDs) provide a visual representation of persistent homology. Each diagram consists of points where the $x$-coordinate indicates the birth time of a feature, that is, when the feature first appears in the filtration, and the $y$-coordinate indicates its death time, that is, when it disappears during the filtration. For more details on TDA and persistent homology, see \cite{TDA}.

The construction of filtrations from signed distance functions and the analysis of the resulting persistence diagrams through functional summaries in \cite{vesnamarcela} enabled to reveal clustering and repulsion tendencies in random set realisations. That framework demonstrated the usefulness of TDA-based techniques in distinguishing between different classes of random set models and in applications such as the analysis of histological images.

This paper continues that line of research. While \cite{vesnamarcela} focused primarily on the design of statistical tests and exploratory tools, here we address the asymptotic properties of persistence diagrams obtained from random set models. Specifically, we establish central limit theorems (CLTs) for functionals of persistence diagrams associated with germ-grain random set models. These results provide the theoretical justification for statistical inference based on persistence diagrams, linking the topological information encoded by persistent homology with probabilistic limit theorems. Our approach builds on recent advances in stabilisation methods for marked point processes that exhibit exponential decay of correlations \cite{prez}, and applies them to marked point processes that generate germ-grain random set models.

The remainder of the paper is organised as follows. Section \ref{sec1} reviews existing theory and necessary background on point processes, marked point processes, and germ-grain random set models. Section \ref{sec2} introduces topological data analysis for random sets and the resulting persistence diagrams, and formalises the notion of $M$-bounded features. Section \ref{sec3} states and proves our main result, namely the central limit theorem for functionals of $M$ bounded persistence diagrams derived from germ-grain models, relying on results from \cite{prez}. Section \ref{sec4} presents statistical tests of goodness of fit for random set models, based on the results from Section \ref{sec3}. Firstly, we do it based on rectangular partition of PD, and then using different summary functions, such as the accumulated persistence function \cite{APF} and the support function of the lift zonoid \cite{vesnamarcela}. We justify them with a simulation study. Section \ref{application} presents the results of applying the method to real-world data, such as histological images of breast tissue. We conclude by discussing the implications of these theorems for statistical applications of TDA to random sets.

\section{Existing theory}\label{sec1}
Before presenting our results, we outline the existing theory required for our work. These definitions are primarily from \cite{chiu:2013,yukich} and \cite{prez}.

Let $\mathcal{P} \subset \R^d$ denote a simple point process. A \textbf{simple point process} is a random element taking values in $\mathcal{N}$, the space of locally finite simple point sets (where no two points coincide) in $\R^d$, equipped with the smallest $\sigma$-algebra such that the number of points in any given Borel set is measurable.

For convenience, we sometimes represent $\mathcal{P}$ as a random set $\mathcal{P} = \{X_i\}$, where $X_i$, $i \geq 1$, are $\R^d$-valued random variables, to emphasise the points of the point process. The fact that it is a simple point process means that $X_i \neq X_j$ for $i \neq j$. Also, note that $\mathcal{P}(B)$ denotes the number of points of the point process that are in the set $B$. We work with stationary point processes, that is, processes that are translation invariant, i.e. $\mathcal{P} + y \overset{d}{=} \mathcal{P}$ for all $y \in \R^d$, where $\mathcal{P} + y$ denotes the translation of $\mathcal{P}$ by the vector $y$.
To state the CLT, we need concepts such as the $k$-point correlation function and reduced Palm expectations. For $(x_1, \dots, x_k) \in (\mathbb{R}^d)^k$, the \textbf{$k$-point correlation function} $\rho^{(k)}$ is determined by the identity $$\mathbb{E}\left(\prod_{1 \leq i \leq k} \mathcal{P}(B_i)\right)=\int \limits_{B_1 \times \dots \times B_k} \rho^{(k)}(x_1, \dots, x_k)\, dx_1 \dots dx_k,$$where $B_1,\dots, B_k$ are mutually disjoint bounded Borel sets in $\mathbb{R}^d$. The $k$-point correlation function provides information about the strength of spatial dependence between points.

The $k$th Palm measure $\mathbb{P}_{x_1,\dots,x_k}$ of $\mathcal{P}$ is the probability distribution of $\mathcal{P}$ conditioned on $\{x_1,\dots,x_k\} \subset \mathcal{P}.$ The \textbf{$k$th Palm expectation} is the expectation corresponding to the Palm probability $\mathbb{P}_{x_1,\dots,x_k}$ and is denoted by $\mathbb{E}_{x_1,\dots,x_k}.$ The reduced Palm versions $\mathbb{P}_{x_1,\dots,x_k}^!$ and their expectation $\mathbb{E}_{x_1,\dots,x_k}^!$ are such that the conditioning points $x_1,\dots,x_k$ are removed from $\mathcal{P}$ and are characterised by $\mathbb{E}^!_{x_1,\dots,x_k}[f(x_1,\dots,x_k;\mathcal{P})]=\mathbb{E}_{x_1,\dots,x_k}[f(x_1,\dots,x_k;\mathcal{P}\setminus\{x_1,\dots,x_k\})],$ for any non-negative measurable function $f$ on $(\R^d)^k \times \mathcal{N.}$ For more details on these concepts, see \cite{yukich}.




We consider point processes that exhibit the so-called exponential decay of correlations. To define this property, we introduce the separation distance between $\mathbf{x} = \{x_1, \dots, x_p\} \subset \R^d$ and $\mathbf{x'} = \{x_{p+1}, \dots, x_{p+q}\} \subset \R^d$ as
\begin{equation}
 \dist(\mathbf{x}, \mathbf{x'}) = \inf_{\substack{i \leq p \\ j \leq q}} |x_i - x_{p+j}|
\end{equation}
where we use the Euclidean distance between points in $\R^d$. For simplicity, we denote $\{x_1, \dots, x_p, x_{p+1}, \dots, x_{p+q}\}$ by $\mathbf{x} \cup \mathbf{x'}$.
\begin{df}
Let $\mathcal{P}$ be a stationary point process in $\R^d,$ such that the k-point correlation function $\rho^{(k)}$ exists for all $k \geq 1.$ Then $\mathcal{P}$ exhibits exponential decay of correlations if there exist $a<1, \phi \colon [0,\infty \rangle \to [0,\infty \rangle $ such that
\begin{itemize}
    \item[1)] $\lim\limits_{t \to \infty} t^n \phi(t)=0$ for all $n \geq 1,$
    \item[2)] $\liminf\limits_{t \to \infty} \log \phi(t)/t^b <0$ for some $b>0,$
    \item[3)] $|\rho^{(p+q)}(\mathbf{x} \cup \mathbf{x'})-\rho^{(p)}(\mathbf{x})\rho^{(q)}(\mathbf{x'})| \leq (p+q)^{a(p+q)}\phi(\dist(\mathbf{x},\mathbf{x'}))$ for any $\mathbf{x}=\{x_1, \dots ,x_p \}, \mathbf{x'}=\{x_{p+1}, \dots ,x_{p+q} \} \subset \R^d.$
\end{itemize}
\end{df}
Intuitively, the fact that the point process exhibits exponential decay of correlations means that the dependence between points that are far apart decreases exponentially with distance, so we can consider them approximately independent. Broadly speaking, such point processes display asymptotic independence over large distances. Examples of point processes with exponential decay of correlations, such as the Mat\'ern cluster process, the determinantal point process, and others, are given in \cite{yukich} and \cite{tgfpp}, together with the proofs. \\
Further on, we need the concept of a marked point process, as we deal with germ-grain random set models, which can be regarded as marked point processes.
\begin{df}
    A marked point process $\tilde{\mathcal{P}}$ is a random sequence $\{(X_i,M_i)\}\subset \R^d\times \mathbb{M}$, where $\mathcal{P}=\{X_i\}$ is a point process in $\R^d$ and the $M_i$ are the marks corresponding to the respective
points $X_i$. The marks $M_i$ belong to a given space of the marks $\mathbb{M}$ which is assumed to be a Polish space. 
\end{df}
\begin{nap}
Equivalently, we can view $\tilde{\mathcal{P}}$ as a random measure $\tilde{\mathcal{P}} = \sum\limits_{x \in \mathcal{P}} \delta_{(x, M(x))}$, meaning that each point of the point process $\mathcal{P}$ is marked with $M(x) \in \mathbb{M}$.
\end{nap}

The concept of exponential decay of correlations in a point process can be generalised to marked point processes as follows:

 \begin{df}
     Let $\Tilde{\mathcal{P}}=\{(X_i,M_i)\}\subset \R^d\times \mathbb{M}$ be a marked point process. We say that it exhibits exponential decay of correlations if there exist $a<1, \phi \colon [0,\infty \rangle \to [0,\infty \rangle $ such that
\begin{itemize}
    \item[1)] $\lim\limits_{t \to \infty} t^n \phi(t)=0$ for all $n \geq 1,$
    \item[2)] $\liminf\limits_{t \to \infty} \log \phi(t)/t^b <0$ for some $b>0,$
    \item[3)] \begin{align*}
   & |\mathbb{E}_{\mathbf{x}\cup \mathbf{x}'}[f(M_1,\dots,M_p)g(M_{p+1},\dots,M_{p+q})]\rho^{(p+q)}(\mathbf{x}\cup \mathbf{x}') \\
    &- \mathbb{E}_{\mathbf{x}}[f(M_1,\dots,M_p)]\rho^{(p)}(\mathbf{x}) \times \mathbb{E}_{\mathbf{x'}}[g(M_{p+1},\dots,M_{p+q})]\rho^{(q)}(\mathbf{x'})|\\
    &\leq (p+q)^{a(p+q)}\phi(\dist(\mathbf{x},\mathbf{x'}))
    \end{align*} for any $\mathbf{x}=\{x_1, \dots ,x_p \}, \mathbf{x'}=\{x_{p+1}, \dots ,x_{p+q} \} \subset \R^d$ and for $||f||_\infty,||g||_\infty\leq 1.$
     
\end{itemize} 
\end{df}
 Note that if the underlying point process exhibits exponential decay of correlations, and the marking is done independently, the marked point process also exhibits exponential decay of correlations. It is also stated in \cite{prez} that this holds if the marking is locally dependent within some deterministic distance.
 
In this work, we present a central limit theorem derived from the persistent homology of a special case of random sets, the so called germ-grain models. For that purpose, we first define the random set.

Let $\mathcal F$ denote the family of closed sets in $\mathbb{R}^d$ and $\mathcal C$ the family of compact sets in $\mathbb R^d.$
\begin{df}
\label{def:randomset}
 Let $(\Omega, \Sigma, P)$ be a probability space.
A mapping $\mathbb{X} :\Omega \to \mathcal{F}$ is a \emph{random closed set} if for every compact set $K \in \mathcal{C}$
$$\left\{ \omega \in \Omega : \mathbb X (\omega) \cap K \neq \emptyset \right\} \in \Sigma.$$
\end{df}

The germ-grain model starts with a marked point process: $$\Tilde{\mathcal{P}} = \{(X_i, M_i)\} \subset \R^d \times \mathbb{M},$$ where the marks are random compact subsets of $\R^d$. The germ-grain model $X$ is derived from this marked point process as the union $X = \bigcup\limits_{n=1}^\infty (X_i + M_i)$.
The points $X_i$ of a point process $\mathcal{P}$ are called germs, and the marks from $\mathbb{M}$ are called grains.
For more details on marked point processes and germ-grain random set models, see \cite{chiu:2013}. In our simulation study, we mainly focus on germ-grain models in which the marks are balls with random radii centred at the germs. Any point process can be used for the underlying point process, but the most commonly applied germ-grain model is the Boolean model, where the underlying point process is a Poisson point process. We also present results for germ-grain models in $\R^2$ other than the Boolean model, but those underlying point processes must satisfy certain conditions, such as exponential decay of correlations.

\section{Topological data analysis for random sets} \label{sec2}

\begin{figure}[ht]
    \centering
    \includegraphics[scale=0.7]{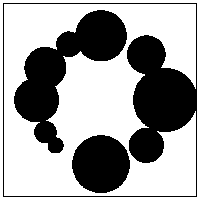}
    \includegraphics[scale=0.7]{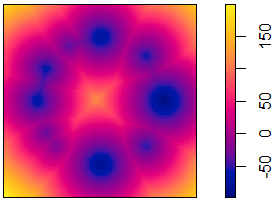} \\
    \caption{Figure showing an example of the set with the heat map of its signed distance function $f_d$. The set is obtained as the union of discs whose centers lie on a circle around the origin with radius $200$. The radii of the discs are $100, 90, 80, 70, 65, 60, 55, 40, 35, 25$.}
    \label{kugle}
\end{figure}

\begin{figure}[ht]
      \centering
       \subfloat[$r=-90$]{\includegraphics[width=0.2\linewidth]{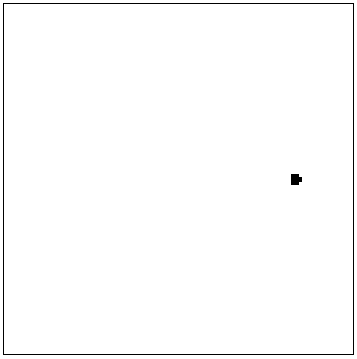}}
      \subfloat[$r=-80$]{\includegraphics[width=0.2\linewidth]{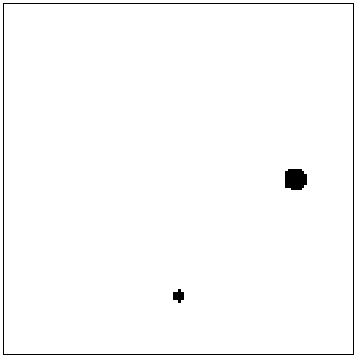}}
       \subfloat[$r=-60$]{
\includegraphics[width=0.2\linewidth]{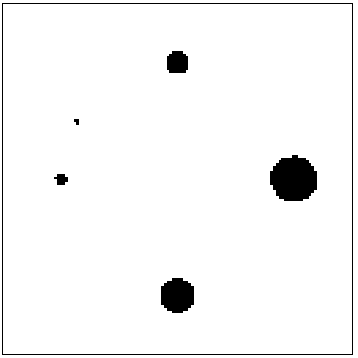}}
       \subfloat[$r=-40$]{
\includegraphics[width=0.2\linewidth]{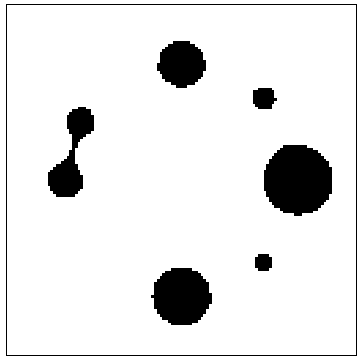}}\\
       \subfloat[$r=-20$]{\includegraphics[width=0.2\linewidth]{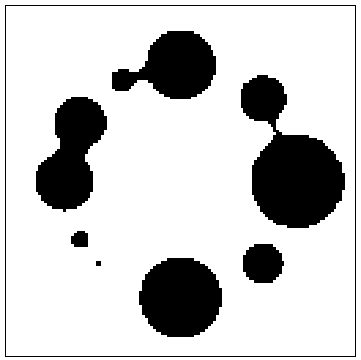}}
       \subfloat[$r=-10$]{\includegraphics[width=0.2\linewidth]{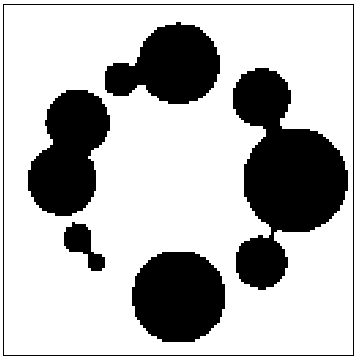}}
       \subfloat[$r=0$]{\includegraphics[width=0.2\linewidth]{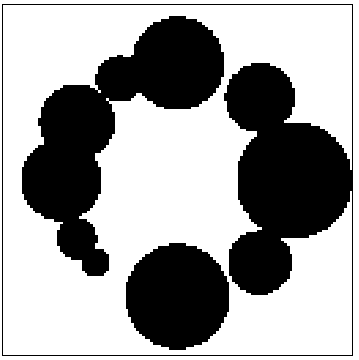}}
       \subfloat[$r=10$]{\includegraphics[width=0.2\linewidth]{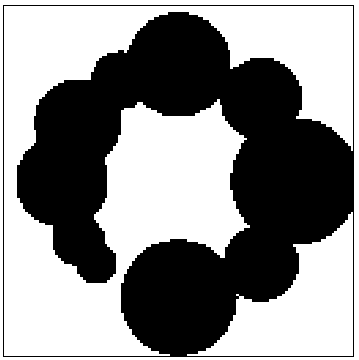}}\\
      \subfloat[$r=30$]{ \includegraphics[width=0.2\linewidth]{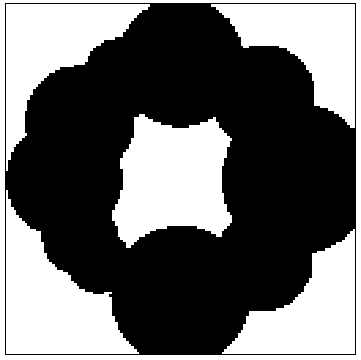}}
      \subfloat[$r=60$]{ \includegraphics[width=0.2\linewidth]{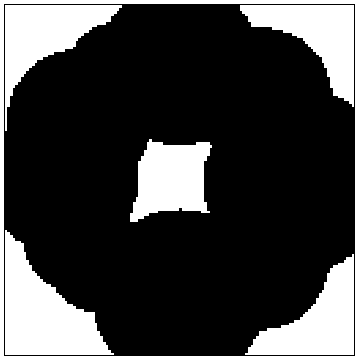}}
             \subfloat[$r=90$]{\includegraphics[width=0.2\linewidth]{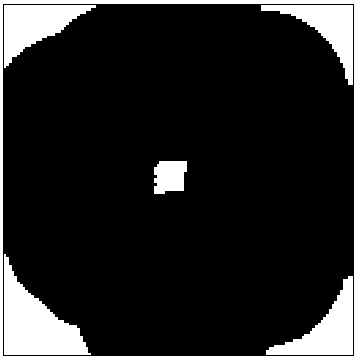}}
      \subfloat[$r=120$]{ \includegraphics[width=0.2\linewidth]{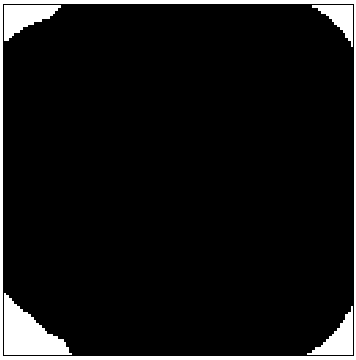}}

      \caption{Sublevel sets of $f_d$ for $r=-90,r=-80,r=-60,r=-40,r=-20,r=-10,r=0,r=10,4=30,r=60,r=90,r=120.$ Note that for $r=0$ we obtain the original set from Figure \ref{kugle}.}
      \label{sublevel}
  \end{figure}

To state our result, we introduce the concept of a persistence diagram from topological data analysis. The core of any TDA method is the construction of a filtration, a parameterised family of sets that grow with the parameter in terms of inclusion. Persistent homology is then derived with respect to this filtration.
We obtain the filtration from the realisation of a random set as follows:
Suppose we work within the observation window $W \subseteq \mathbb{R}^2$ and observe the realisations of our random set.
Let us consider the signed distance function to a set $S$, defined as  $f_d: W \to \mathbb R,$
\[f_d(x)=\left\{\begin{array}{ll}
d(x,S), & x \notin S,\\
-d(x,W\diagdown S) & x \in S,
\end{array}\right. \]
where $d(x,S)$ is the distance from a point $x$ to a set $S$ calculated using the Euclidean distance. If we consider $S_r,$  the sublevel sets of $f_d,$ i.e.
$$S_r=f_d^{-1}\left(\left\langle-\infty,r \right]\right),$$
we obtain a non-decreasing filtration  $\{S_r\}_{r \in \mathbb{R}}$, from which we can construct a persistence diagram. Figure \ref{kugle} shows an example of a set together with the heat map of its signed distance function $f_d.$ We can imagine that as the parameter changes, the set grows through filtration starting from an empty set. Figure \ref{sublevel} illustrates the sublevel sets of the signed distance function from Figure \ref{kugle}, for $r=-90,r=-80,r=-60,r=-40,r=-20,r=-10,r=0,r=10,4=30,r=60,r=90,r=120.$ Note that for $r=0$ we obtain the original set.

 As it is not crucial for this work, we proceed without providing a formal definition of homology and persistent homology (see \cite{TDA}), but we explain what the persistence diagram can capture about the random set. Persistent homology is used to track the persistence of topological features in the filtration as the parameter changes. By topological features, we mean the connected components of the objects and the loops formed by the objects. As the parameter changes, different features in the filtration appear and disappear. When a connected component appears for the first time at the sublevel $S_r$ of the filtration, we say that it is born at time $r$. If a connected component connects to the other connected component at the sublevel $S_r$, we say that the feature dies at time $r$. When two components merge, the younger one is considered to die, while the older one continues to live until it merges with an even older component. The oldest component dies when all components have merged into a single component. For loops, we say that a loop is born at time $r$ if it appears in the filtration at parameter $r$. A loop can only be formed when two connected components merge. The death time of a loop is the first time the hole is completely filled. See Figure \ref{sublevel} for a clearer visualisation of the appearance and disappearance of features. We consider connected components as $0$-dimensional topological features and loops as $1$-dimensional topological features. Features that persist over time are usually significant when describing the dataset.

\begin{figure}[H]
    \centering
    \includegraphics[width=4cm, height=4cm]{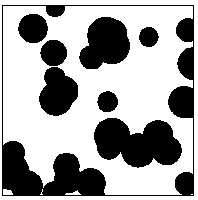} 
    \includegraphics[width=4.8cm, height=3.95cm]{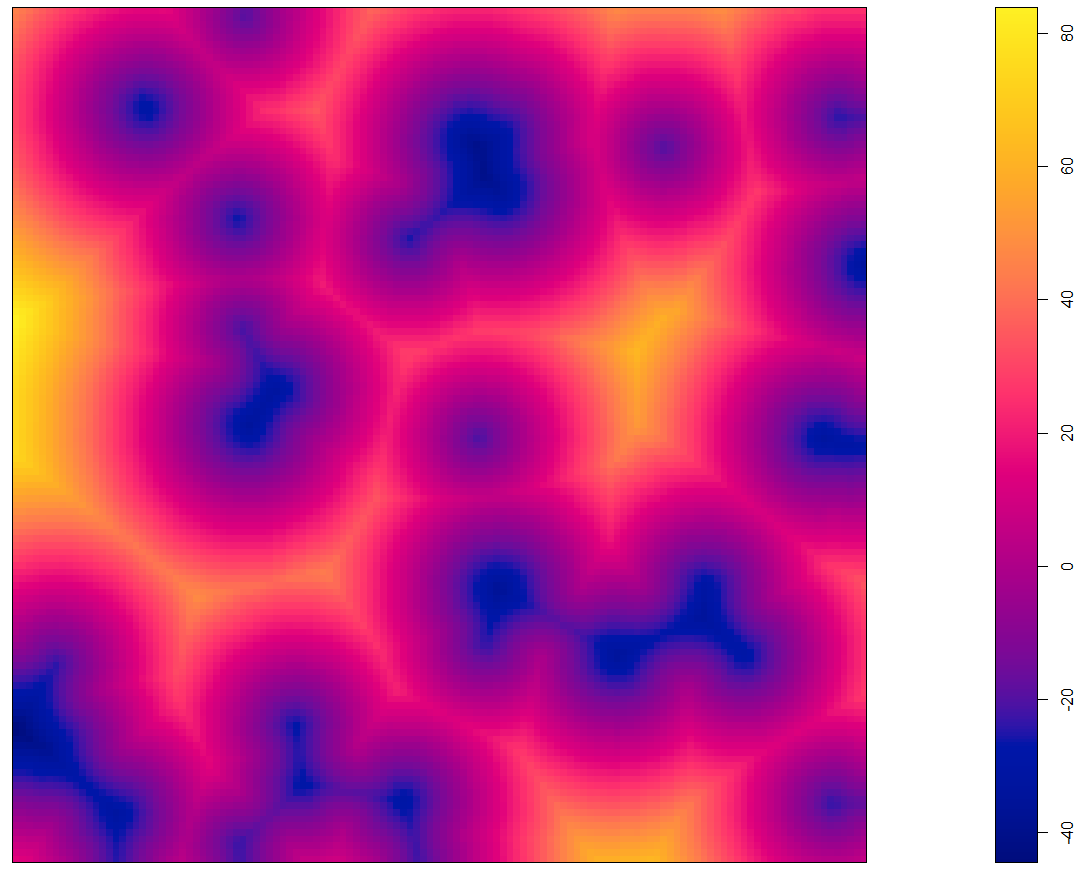}\\
    \includegraphics[width=4.5cm, height=4.2cm]{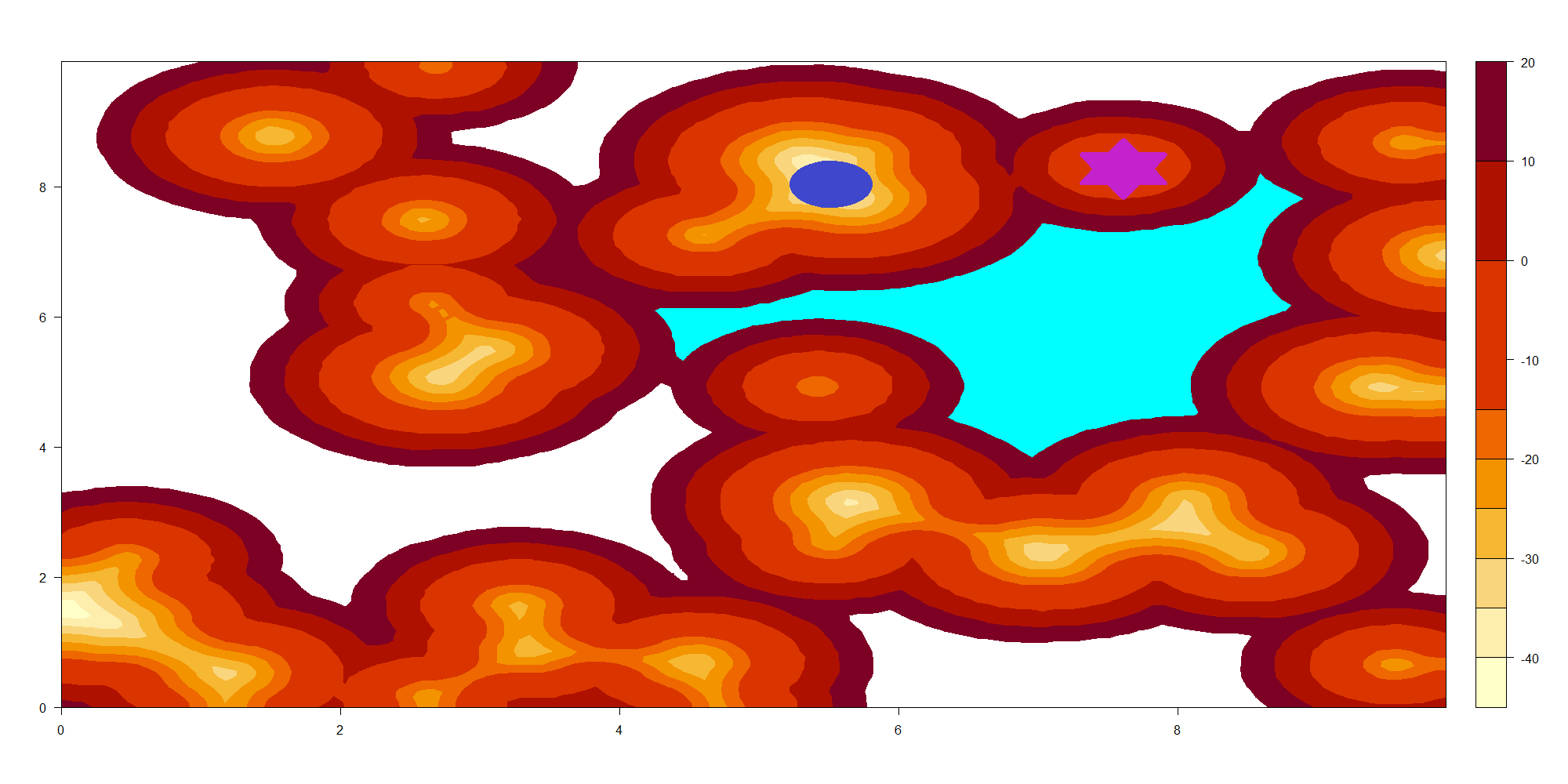}
    \includegraphics[width=4 cm, height=3.95 cm]{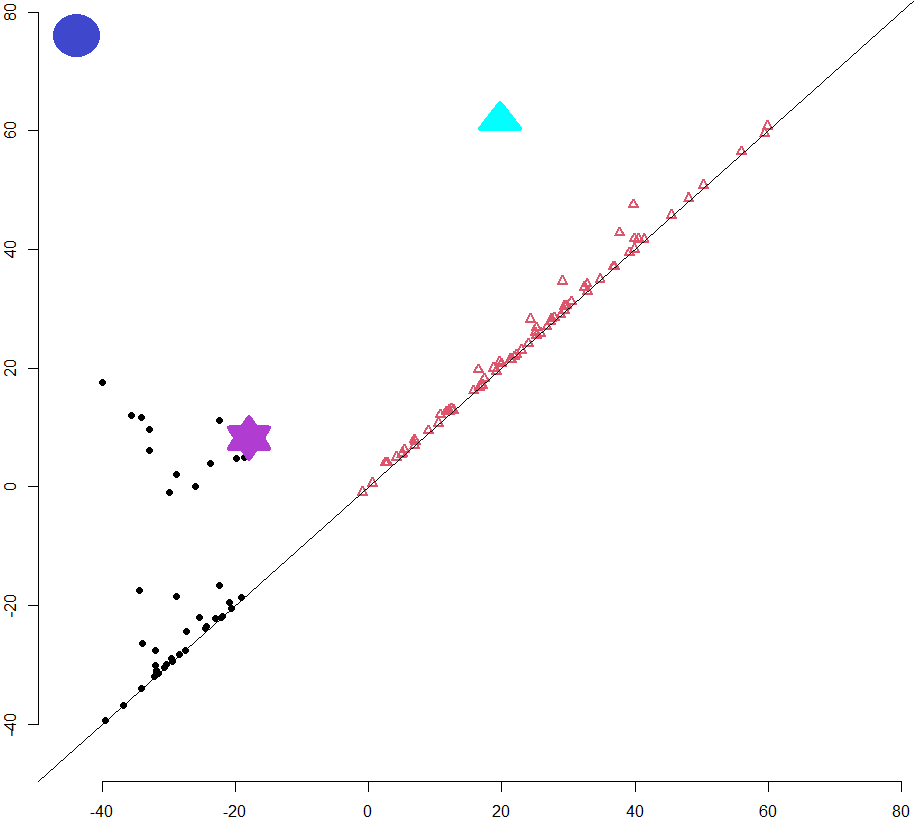}
    \caption{Top left: Realisation of the germ-grain model of the random set. Top right: Signed distance function of the realisation of the random set shown in the top-left figure. Bottom left: Sublevel set of the signed distance function for the thresholds $-45, -40, -35, -30, -25, -20, -20, 1, 10,$ and $20$. Bottom right: PD corresponding to the sublevel sets of the signed distance function (black dots, purple star, and blue dot represent the 0-dimensional homology; pink triangles and the blue triangle correspond to the 1-dimensional homology).}
    \label{fig:uvod}
\end{figure}

 For a fixed dimension $q \in \{0,1\}$, the persistence diagram $PD^q$ graphically represents the points with abscissae equal to the birth times and ordinates equal to the death times for all $q$-dimensional features of a random set $Z$. Note that some features may be born at the same time (they appear at the same sublevel) and may die at the same time (they disappear at the same sublevel). Without going into the theoretical details of constructing persistence diagrams, we illustrate it on a simple example in Figure \ref{fig:uvod}.
The first row of Figure \ref{fig:uvod} shows the realisation of the germ-grain model (left) and its corresponding signed distance function (right). In the heat map of the signed distance function $f_d$, the local minima appear near the centers of the grains, that is, close to the points of the underlying point process $\mathcal{P}$. If the grains do not overlap, the local minima occur precisely at the points of the underlying point process, with values equal to the negatives of the corresponding ball radii. The local maxima of the signed distance function appear near the centers of the holes enclosed by the components. We denote the birth times by $b_i$ and the death times by $d_i$. Note that $b_i$ and $d_i$ correspond to the values of the signed distance function at its stationary points, that is, at its local or global maxima and minima. The left image in the second row of Figure \ref{fig:uvod} shows sublevel sets of the signed distance functions for the thresholds $-45, -40, -35, -30, -25, -20, 1, 10,$ and $20$. The blue dot represents the global minimum of the signed distance function, which is the first-born component.

The bottom-right image in Figure \ref{fig:uvod} shows the persistence diagram obtained from the filtration of the signed distance function of the random set realisation depicted in the same figure. The 0-dimensional persistent homology class is indicated by black dots, with a purple star marking the birth and death coordinates of the component highlighted by the purple star in the bottom-left image, and a blue dot representing the oldest connected component. Triangles denote the 1-dimensional persistent homology class, with the blue triangle indicating the birth and death coordinates of the blue hole shown in the bottom-left image.

The component marked by the purple star is born at $-20$ and dies at sublevel $10$, when it merges with the component that began to grow at the blue dot. Additionally, the light blue hole appears at sublevel $20$ and disappears at sublevel $60$. 

\begin{figure}[h!]
 \centering
 \includegraphics[scale=0.5]{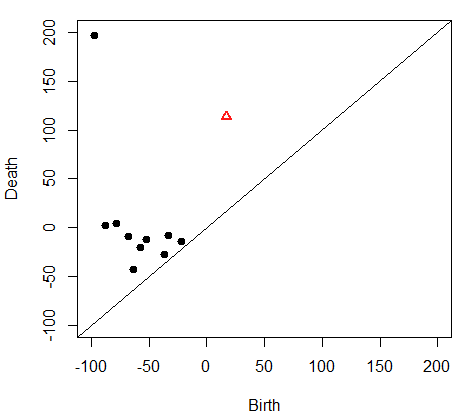}
 \caption{PD of the set from Example \ref{kugle}.}
 \label{PD_kugle}
\end{figure}

Figure \ref{PD_kugle} shows the PD of the set from Example \ref{kugle}. Notice that the number of black points corresponds to the number of discs, with the first coordinate being the negative value of their radius. The red triangle represents the hole encircled by the balls.

The persistence diagram can be regarded as an empirical measure
\[
PD^{q}(Z)=\sum\limits_{i \in \mathcal I_q }\delta_{(b_i,d_i)},\] where $\delta$ stands for the Dirac delta
measure which equals $1$ if $(b_i,d_i)$ is in the persistence diagram and $0$ otherwise.  $\mathcal I_q$ is an index set over all $q$-dimensional features and $b_i,$ $d_i$ are the birth and death times of the $i$-th feature, with $q \in \{0,1\}$. Since some features may have the same birth and death times, the point $(b_i,d_i)$ may have a multiplicity $c_i$ greater than 1, which is not observable in the plots of PDs.

Henceforth, for practical reasons, we observe the signed distance function corestricted to $[-r_f,r_f]$ where $r_f$ is a large fixed deterministic radius. Consequently, all births and deaths are within the segment $[-r_f,r_f]$ so the persistence diagram is on the given segment.

Also, we do not consider arbitrarily large topological features. Spatially bounded features are sufficiently useful for simulated point processes and for the real datasets we use. We do this because some localisation conditions are required for the proof of the CLT, which will be discussed in Section \ref{CLTmpp}. Therefore, we introduce the concept of $M$-bounded features, following the intuition of \cite{tgfpp}. Regarding the birth, $M$-bounded $0$-dimensional features are those connected components that are born on a sublevel set $S_r$ of a signed distance function with $r>-M.$ The $0$-dimensional features die on a sublevel set of the signed distance function where the set of germs contained in the feature has diameter greater than $M$, or when the component merges with an older component. Concerning the holes, if the diameter of a set of germs whose marks created it is larger than $M,$ we do not consider it to be born until, at some sublevel, it has a diameter smaller than $M$. In this way, we ensure that only points that are not distant from each other have influence on birth and death of associated features, which will secure the aforementioned localisation.


We denote the persistence diagram that includes only $M$-bounded features by $$PD^{M,q}(Z)=\sum\limits_{i \in \mathcal I_{M,q} } \delta_{(b_i^M,d_i^M)},$$ where the $b_i^M$ and $d_i^M$ are births and deaths of $M$-bounded features respectively.

\begin{nap}
    Note that methods from topological data analysis are not restricted to dimension $2$. The same concepts can be defined for any dimension $d\geq2$, but the interpretations of the features become more complex. Apart from dimension $2$, in dimension $3$, topological features remain interpretable: $0$-dimensional features are connected components, $1$-dimensional features are tunnels, and $2$-dimensional features are voids or cavities, that is, enclosed empty regions. We do not consider higher dimensions in this work. In the following sections, we present results for dimension $2$ and assume that they can be generalized to dimension $3$ in a similar manner.
\end{nap}

\section{Central limit theorem for germ-grain random set models} \label{sec3}

Our main goal is to prove the central limit theorem for functionals of persistence diagrams in germ-grain random set models with respect to a growing observation window. We rely on the frameworks of \cite{yukich}, \cite{tgfpp}, and \cite{prez}, and apply these results in our context.

Assume we have the germ-grain model $\tilde{\mathcal{P}}$, with the underlying point process $\mathcal{P}$ in $\R^2$ and marks that are random sets in $\R^2$, which are strictly convex and compact. A set is strictly convex if it is convex and every line segment connecting two distinct boundary points lies entirely in the interior. Those types of sets have a unique largest inscribed ball (inball), see \cite{inball}. We assume the marking is done so that the center of the inball coincides with the point from the underlying point process. Also, we limit the radius of those inballs to be less than some large deterministic $R_{max}$. Assume that the following relation holds: $R_{max}<M<r_f.$\\
Let $\mathcal{P}_n=\mathcal{P} \cap W_n$, where $W_n=\left[-\frac{\sqrt{n}}{2}, \frac{\sqrt{n}}{2}\right]^2$. We now observe only the points of $\mathcal{P}_n$ together with their marks, and denote this random set by $\tilde{\mathcal{P}}_n$. \\
To state the desired CLT precisely, we let $$\langle f, PD^{M,q}(\tilde{\mathcal{P}}_n)\rangle = \int\limits_{[-r_f, r_f]^2} f(b, d)\, PD^{M,q}(\tilde{\mathcal{P}}_n)(\mathsf{d}b,\mathsf{d}d)=\sum_{i \in \mathcal I_{M,q}(\tilde{\mathcal{P}}_n)} f(b_i^M, d_i^M)$$ denote the integral of a bounded measurable function $f\colon [-r_f, r_f]^2 \to \R$ with respect to the measure $PD^{M,q}(\tilde{\mathcal{P}}_n).$
\subsection{General CLT for marked point processes}\label{CLTmpp}

Following the framework of \cite{prez}, we assign a score to each point in the point process via the score function $\xi$. To apply those results, the score should depend only on the neighbouring points and marks; that is, there is only "local influence" on the value of the score. We assign a real-valued score, denoted by $\xi(\tilde{x},\tilde{\mathcal{P}})$, to each point $\tilde{x}$ from $\tilde{\mathcal{P}}$. Note that those results are more general than what we use here. In our case, we consider score functions with a bounded radius of stabilisation. Further on, $B_r(x)$ denotes the ball of radius $r$ centered at $x.$

\begin{df}
    Given a score function $\xi$, locally finite marked point process $\tilde{\mathcal{X}}$ and $\tilde{x}\in \tilde{\mathcal{X}},$ define the radius of stabilization $ R^{\xi}(\tilde{x},\tilde{\mathcal X})$ as the smallest $r \in \Z_{\geq 0}$ such that $$\xi(\tilde{x},\tilde{\mathcal{X}}\cap B_r(x))=\xi(\tilde{x},(\tilde{\mathcal{X}}\cap B_r(x))\cup (\tilde{\mathcal{A}}\cap B_r^c(x)))$$ for all locally finite $\tilde{\mathcal A} \in \R^2 \times \mathbb M.$ If no such finite $r$ exists, we set $ R^{\xi}(\tilde{x},\tilde{\mathcal{X}})=\infty.$
\end{df}

We say that the radius of stabilisation is bounded if $ R^{\xi}(\tilde{x},\tilde{\mathcal X})$ is uniformly bounded over all $\tilde{\mathcal{X}}$, and all $\tilde{x}\in \tilde{\mathcal{X}}.$
The fact that the radius of stabilisation is bounded means that points far from $\tilde{x}$ do not influence the score. If this holds, then the "BL-stabilisation" from \cite{prez} follows, as it requires that the influence of points far from $\tilde{x}$ decreases with distance.

The theorem we rely on is stated more generally in \cite{prez}, but for our purposes, we present it only for bounded radii of stabilisation.

\begin{tm}\label{CLTsf}
    Let $\tilde{\mathcal{P}}$ be a marked point process that exhibits exponential decay of correlations, and let $\xi$ be a score function that has bounded radius of stabilisation. Let $\xi(\tilde{x_1},\tilde{\mathcal{P}_n})$ satisfy the $p$-moment condition, that is $$\sup\limits_{1\leq n <\infty}\sup\limits_{1\leq q \leq p} \sup_{x_1,\dots,x_q \in W_n}\mathbb E_{x_1,\dots,x_q}[|\xi(\tilde{x_1},\tilde{\mathcal{P}_n})|^p]<\infty.$$ For $H_n:=\sum\limits_{x \in \mathcal{P}_n}\xi(\tilde{x},\tilde{\mathcal{P}_n}),$ if $\Var(H_n)=\Omega(n^\nu),\nu>0,$ then $$\frac{H_n-\mathbb{E}(H_n)}{\sqrt{\Var(H_n)}}\overset{D}{\longrightarrow}N(0,1).$$
\end{tm}

\begin{nap}
    The notion  $X=\Omega(n)$  means that $X$ grows at least as fast as $n,$ i.e. $\operatorname{liminf}\limits_{n\to \infty}\frac{X}{n}>0.$ 
\end{nap}

Examples of processes that satisfy the $p$-moment condition can be found in \cite{yukich}, for instance, when the underlying point process is determinantal or permanental. It is also stated in \cite{prez} when the variance condition is satisfied.

\subsection{Score functions and CLT in our context}
We use the previous theorem to prove the central limit theorem for the integral of a bounded measurable function with respect to the measure $PD^{M,q}(\tilde{\mathcal{P}}_n)$. Inspired by the work in \cite{tgfpp}, our aim is to express  the quantity $$\langle f, PD^{M,q}(\tilde{\mathcal{P}}_n)\rangle=\sum\limits_{i\in \mathcal I_{M,q}(\tilde{\mathcal{P}}_n)}f(b_i^M,d_i^M)$$ in the form $$\sum\limits_{z \in \tilde{\mathcal{P}}_n}\xi(\tilde{z},\tilde{\mathcal{P}}_n)$$ for a suitable score function $\xi(\tilde{z}, \tilde{\mathcal{P}}_n)$ so that we can apply Theorem \ref{CLTsf}. We need to transform indexing over features into indexing over points from the marked point process. This means that for each point in PD, we must determine which point in the marked point process produced it, that is, which point is responsible for its occurrence (for the birth or death of a feature).
We say that the $(\star)$ condition is satisfied if for every $p\geq 1,$ $$\sup\limits_{1 \leq q \leq p} \sup\limits_{x_1,\dots,x_q \in \R^2} \mathbb{E}^!_{x_1,\dots,x_q}[(\tilde{\mathcal{P}}(W_1))^p]<\infty. $$ 
\begin{tm}\label{CLTPD}
   Let $M>0, q\in \{0,1\}$ and $f\colon [-r_f,r_f]^2 \to \R$ be a bounded measurable function and assume that $\tilde{\mathcal{P}}$ exhibits an exponential decay of correlations and satisfies the $(\star)$ condition. If $\operatorname{Var}(\langle f, PD^{M,q}(\tilde{\mathcal{P}_n})\rangle)=\Omega(n^\nu)$ for $\nu>0,$ then  $$\frac{\langle f, PD^{M,q}(\tilde{\mathcal{P}_n})\rangle-\mathbb{E}[\langle f, PD^{M,q}(\tilde{\mathcal{P}_n})\rangle]}{\sqrt{\Var(\langle f, PD^{M,q}(\tilde{\mathcal{P}_n})\rangle)}}$$ converges in distribution to a standard normal random variable as $n \to \infty.$
\end{tm}

Firstly, we do it for $q=0.$ \\We can express $\sum\limits_{i\in \mathcal I_{M,0}(\tilde{\mathcal{P}_n})}f(b_i^M,d_i^M)$ in the form $\sum\limits_{z \in \mathcal{P}_n}\xi_0(\tilde{z},\tilde{\mathcal{P}_n})$ where we put $\xi_0(\tilde{z},\tilde{\mathcal{P}_n})=\sum\limits_{i\in \mathcal I_{M,0}(\tilde{\mathcal{P}_n})} \mathbbm{1}_{\{ \tilde z \text{ gives birth to feature i} \}}  f(b_i^M,d_i^M).$ \\ We need to clarify what it means for the point $\tilde z$ from $\tilde{\mathcal{P}_n}$ to be responsible for the occurrence of a point $(b_i^M, d_i^M)$ in the persistence diagram (i.e. when we consider that it gives birth to feature $i$).

\begin{itemize}
    \item We say that $\tilde{z}$ gives birth to feature $i$ if there is a local minimum of the signed distance function, with value $b_i^M$, in the inball of the mark attached to a point $z$ from $\mathcal{P}_n$. If the minimum lies in more than one inball, we assign responsibility for the birth of feature $i$ to the point whose mark has the largest inball radius among those that have not yet been assigned responsibility for any other feature.
    
    
    \item In the simplest case where there is no overlapping of grains, it means that the signed distance function has a local minimum, valued $b_i^M,$ at the center of an inball of the grain. The value $b_i^M$ is the negative radius of that inball, and the points where the local minima occur coincide with the germs.
    
\end{itemize} 

We need to check whether our score function satisfies the conditions of Theorem \ref{CLTsf}. The key observation is that the point $\tilde z$ in our marked point process can be responsible for at most one point in $PD$, meaning that one point in the marked point process can give birth to at most one feature.


    The $p$-moment condition holds since we have: $$|\xi_0(\tilde{z},\tilde{\mathcal{P}_n})| \leq ||f||_{\infty}<\infty$$ which implies that $$\mathbb{E}_{x_1,\dots,x_q}[|\xi_0(\tilde{z},\tilde{\mathcal{P}_n})|^p] \leq \mathbb{E}_{x_1,\dots,x_q}||f||_{\infty}^p=||f||_{\infty}^p < \infty$$ for all $p>1,$ and $1\leq q\leq p.$ \\ The radius of stabilization is bounded since we have
    \begin{equation}
    \begin{aligned}
        \sum_{i\in \mathcal I_{M,0}(\tilde{\mathcal{P}}\cap B_r(z))} \mathbbm{1}_{\{z \text{ gives birth to feature i} \}} f(b_i^M,d_i^M)=\\ \sum_{i\in \mathcal I_{M,0}(\tilde{\mathcal{P}}\cap B_r(z)) \cup (\tilde{\mathcal A} \cap B_r^c(z))} \mathbbm{1}_{\{z \text{ gives birth to feature i} \}} f(b_i^M,d_i^M)=f(b_j^M,d_j^M)
        \end{aligned}
    \end{equation} for some $j$ and all $\tilde{\mathcal A}$ that are locally bounded. This equality holds for all $r > M$. The main reason for these equalities is that one point of the marked point process can be responsible for at most one point in $PD$, as mentioned above, and because we have reduced to $M$-bounded features. \\


    Now, we can proceed in the same way for $q=1.$ \\
    We can express $\sum\limits_{i\in \mathcal I_{M,1}(\tilde{\mathcal{P}_n})}f(b_i^M,d_i^M)$ in the form $\sum\limits_{z \in \mathcal{P}_n}\xi_1(\tilde{z},\tilde{\mathcal{P}_n})$ where we put $\xi_1(\tilde{z},\tilde{\mathcal{P}_n})=\sum\limits_{i\in \mathcal I_{M,1}(\tilde{\mathcal{P}_n})} \mathbbm{1}_{\{z \text{ gives birth to feature i} \}} f(b_i^M,d_i^M).$ \\

    A hole can occur only when, at some level set of the signed distance function, two components meet and connect. The key observation is that one point of a marked point process can give birth to at most $\tilde{\mathcal{ P}}(W_n) - 1$ holes. This holds because each hole is caused by a point of the process, so there can be no more holes than points.\\ If $\tilde z$ is responsible for the birth of a hole at time $r > 0$, then there exists $\tilde{z}'$ in $\tilde{\mathcal{P}_n}$ at distance $2r$ from $\tilde z$. If the birth of a hole occurs at time $r < 0$, it means that the hole is formed between overlapping grains. Since we restricted the values of the signed distance function to $[-r_f, r_f]$, the following holds: 
$$|\xi_1(\tilde{z}, \tilde{\mathcal{P}}_n)| \leq \tilde{\mathcal{P}}_n(B_{4r_f}(z)) \|f\|_{\infty}.$$
    The $p$-th moment condition is satisfied after applying stationarity and the $(\star)$ condition.
    The radius of stabilization is bounded because we reduced to $M$- bounded features.

    \begin{nap}
       Suppose we consider only independently marked point processes, meaning the grains are independent and are also independent of the germs. Then, it is sufficient to assume in Theorem \ref{CLTPD} that the underlying point process $\mathcal{P}$ has exponential decay of correlations, as this directly implies that $\tilde{\mathcal{P}}$ satisfies the same condition due to independent marking.
    \end{nap}

\section{Simulation study}\label{sec4}
In our simulation study, we mostly use germ-grain models in which marking is done independently. We also use more complex models that allow dependence between the germs and grains, such as the Quermass-interaction process.

\begin{df}[\cite{kendall:1999}]
\label{def:quermass}
 Let $\mathbb X$ be the Boolean model with 
$\mathbf B_1$ being a disc in $\mathbb R^2$ with the random radius.
\emph{Quermass-interaction process} is a random set 
whose probability measure 
is absolutely continuous with respect to the 
probability measure of  $\mathbb X$  with density given by
\begin{align}
	f_{\theta}(\mathbf b)
    & = 
    \frac{1}{c_{ \theta}}
    \exp\{
    	\theta_1A(  U_{\mathbf b})
        + \theta_2L(  U_{\mathbf b})
        + \theta_3\chi(  U_{\mathbf b})
    \}
    \label{eq:quermass_den}
\end{align}
for each finite disc configuration
$\mathbf b = \{\mathbf b_1\ldots,\mathbf b_n\}$, where
$A=A(  U_{\mathbf b})$ is the area,
$L=L(  U_{\mathbf b})$ is the perimeter,
$\chi=\chi(  U_{\mathbf b})$ is Euler-Poincar\'{e} characteristic
(i.e. the number of connected components minus the number of holes)
of the union $U_{\mathbf b} = \cup_{i=1}^n \mathbf b_i$, 
$\theta = (\theta_1,\theta_2,\theta_3)$ is 3-dimensional vector 
of parameters and $c_{\theta}$ is the normalising constant.
\end{df}

By adjusting the parameters $\theta_1, \theta_2,$ and $\theta_3$ in the density function given in equation \ref{eq:quermass_den} of the Quermass-interaction process, various interaction patterns between grains can be produced. When $\theta_1$ is positive and $\theta_2$ negative, the resulting configurations tend to have larger areas and smaller perimeters than the corresponding Boolean model, leading to clustering of the grains. Conversely, setting $\theta_1$ negative and $\theta_2$ positive yields configurations with smaller areas and larger perimeters, producing non-overlapping configurations, and thus repulsion of the grains.
It is shown in \cite{gibbs} that models with $\theta_3 = 0$, and where the radii of the balls have support on some positive interval, satisfy exponential decay of correlations under certain assumptions. The process must be in the sparse regime, below the percolation threshold of the corresponding Boolean model; that is, there should be no infinite components, and only local clusters should occur.

Our simulation studies are based on the following random set models, all defined with disc centres within a $25\times25$ window:

\begin{itemize}
    \item Boolean model in which the grains are discs with random radii, the intensity of the germs is $0.4$, and the radii are uniformly distributed over the interval $(0.5, 1)$ (referred to as the \emph{Boolean} model).
    \item Boolean model in which the grains are ellipses, with their axes fixed and parallel to the two Cartesian axes. The major and minor axes are random, with the intensity of the ellipse centres equal to $0.4$, the semi-major axes uniformly distributed on the interval $(0.5, 1)$, and the semi-minor axes uniformly distributed on the interval $(0.2, 0.7)$ (referred to as a \emph{Boolean ellipse} model).
    
    \item Quermass-interaction process (see Definition \ref{def:quermass})  in which we used parameters $\theta_1 = 0.62, \theta_2 = -0.86$ and $\theta_3 = 0.7$ with respect to the above-mentioned Boolean model. The model produces realisations with larger area, smaller perimeter, and higher Euler–Poincar\'{e} characteristic than the Boolean model, thereby favouring the formation of clusters. It is referred to as a \emph{cluster} model.

 \item Quermass-interaction process (see Definition \ref{def:quermass}) with parameters $\theta_1 = -1$, $\theta_2 = 1$, and $\theta_3 = 0$ with respect to the same random-disc Boolean model. The model yields realisations with reduced area and increased perimeter compared to the Boolean model, resulting in typical realisations consisting of smaller, non-overlapping parts. It is referred to as a \emph{repulsive} model.

\item Model in which the germs form a Mat\'ern cluster point process \cite{matern} that we construct using parent points forming a homogeneous Poisson point process that has intensity $0.4$. Each parent point is replaced by a random cluster of offspring points forming a homogeneous Poisson process with intensity $\frac{1}{4\pi}$ within the disc of radius $2$ centred at the parent point. The grains are again discs with uniformly distributed radii on the interval $(0.5, 1)$. They are independent of the process of grains. This model is referred to as a \emph{Matern cluster} later in the text.

\item Model in which the germs form a Baddeley-Silverman cell process \cite{baddeley}. This is a random point process with the same second-order properties as the homogeneous Poisson point process, but with more complex higher-order interactions. The model is obtained by dividing the space into equal rectangular cells, with each cell containing a random number of points that are uniformly and independently distributed. This number is $0$, $1$, or $28$ with probabilities $\frac{1}{28}$, $\frac{26}{27}$, and $\frac{1}{765}$, respectively. Setting the cell size to $\frac{25}{\sqrt{280}}$ ensures an average realisation intensity of $0.4$. The grains are discs with radii uniformly distributed on the interval $(0.5, 1)$, independent of the germ process. This model is referred to later in the text as a \emph{Cell process}.

 \item The model in which germs form a determinantal point process of Bessel type \cite{DPP} with intensity $0.4$, 
where the repulsion between points is governed by a Bessel-type kernel with parameters $\sigma = 0$ and $\alpha = \frac{1}{\sqrt{0.4\pi}} - 0.01$. 
The grains are discs with radii independently drawn from a uniform distribution on $(0.5, 1)$. 
In the following, this model is reffered to as \emph{DPP}.
\end{itemize}

\begin{figure}[h!]
\centering
 \includegraphics[scale=0.4]{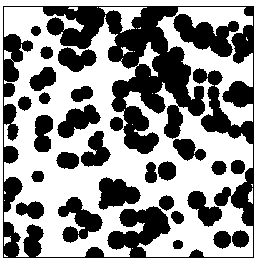}
   \includegraphics[scale=0.4]{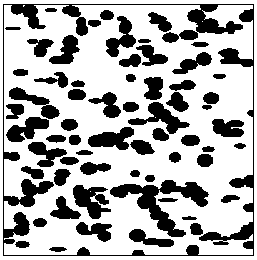}
   \includegraphics[scale=0.4]{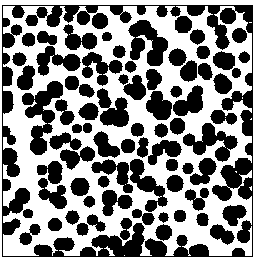}
   \includegraphics[scale=0.4]{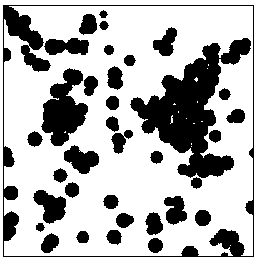} \\
   \includegraphics[scale=0.4]{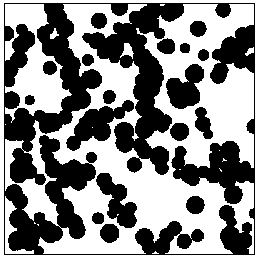}
   \includegraphics[scale=0.3]{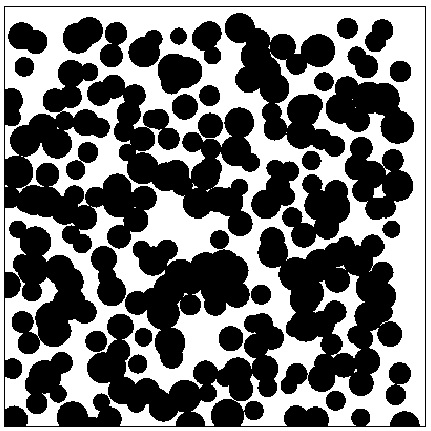}
   \includegraphics[scale=0.3]{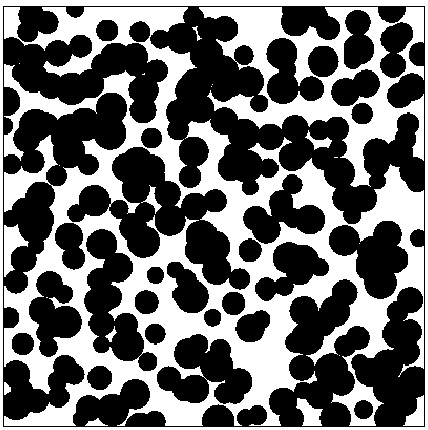 }\\
    \caption{One realisation of each random set process in this order: Boolean, Boolean ellipse, repulsive, cluster, Matern cluster, Cell and DPP.}
       \label{fig:procesi}
\end{figure}

\begin{figure}[h!]
    \centering
   \includegraphics[scale=0.15]{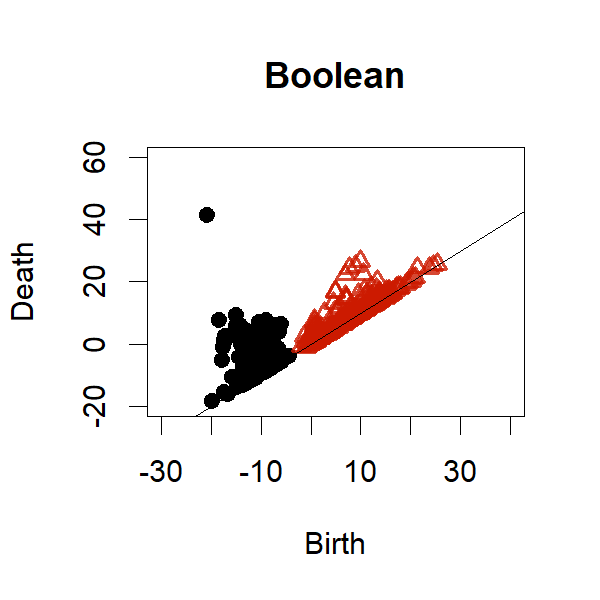}
   \includegraphics[scale=0.15]{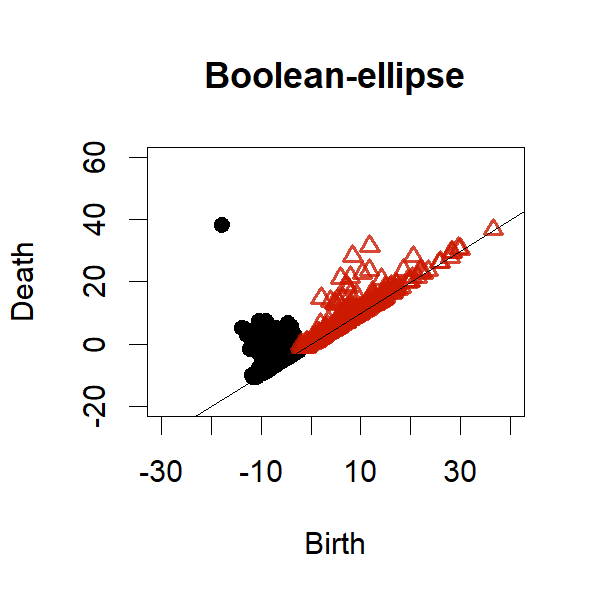}
   \includegraphics[scale=0.15]{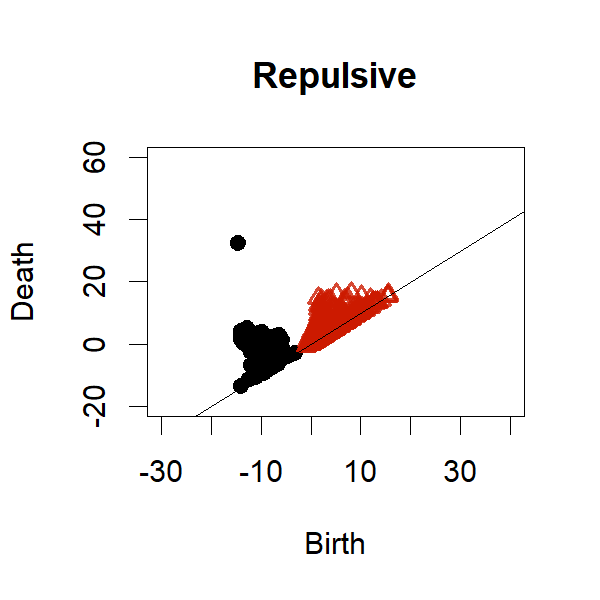}
   \includegraphics[scale=0.15]{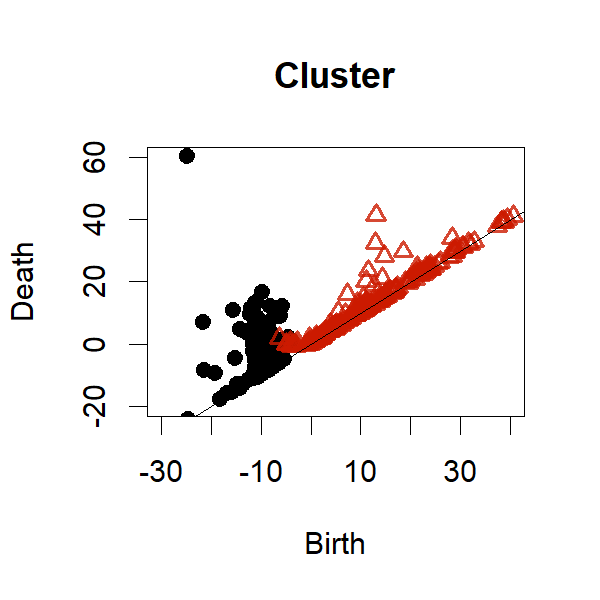}
   \includegraphics[scale=0.15]{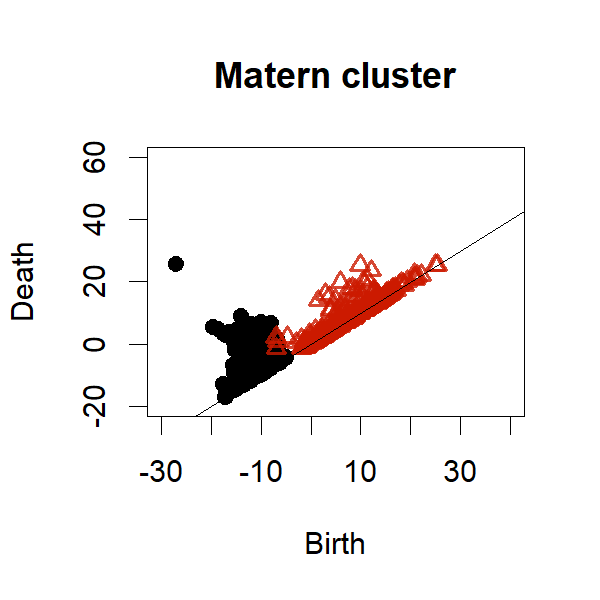}
   \includegraphics[scale=0.15]{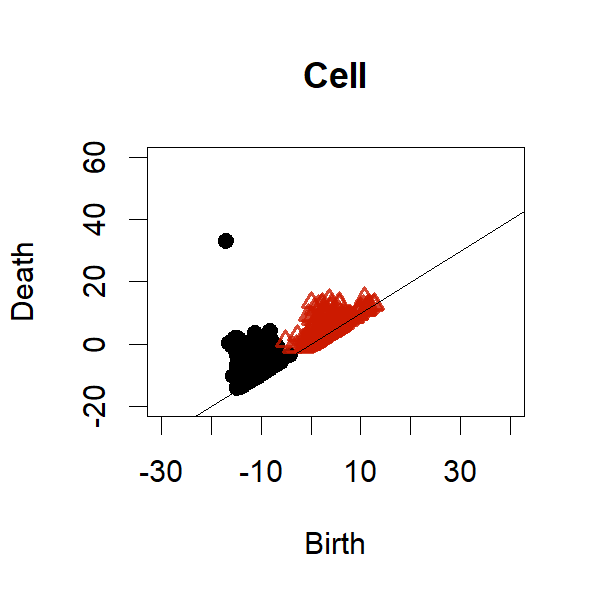}
   \includegraphics[scale=0.15]{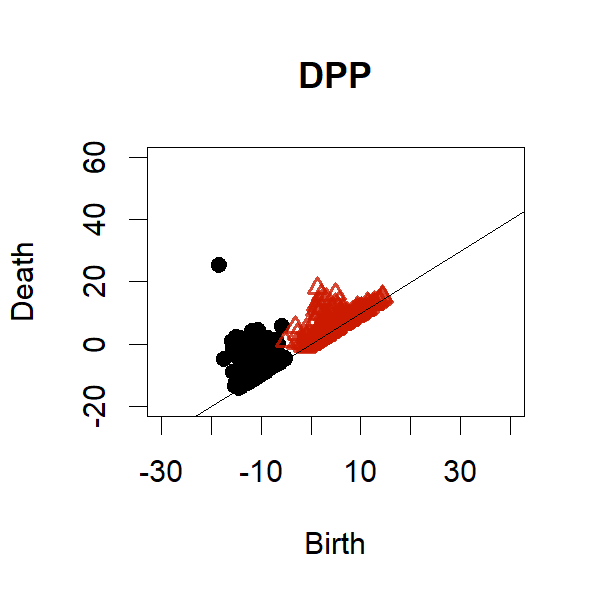}
    \caption{Persistence diagram of realisations used in simulation study given in this order: Boolean, Boolean ellipse, Repulsive, Cluster, Matern cluster, Cell, and DPP.}
    \label{PD}
\end{figure}

A realisation of each of the mentioned random set models is shown in Figure \ref{fig:procesi}. Figure \ref{PD} presents the persistence diagrams of these realisations. These images are taken from \cite{vesnamarcela}.

We used the library \verb{spatstat{ \cite{spatstat} in R for the simulation study, except for the simulation of the Quermass-interaction process, where we used the algorithm from \cite{moeller:2008}. We computed PDs using the function \verb{gridDiag{ from the package \verb{TDA{ \cite{TDA_R}.

To apply Theorem \ref{CLTPD} to the models listed above, one should verify whether the assumptions of the theorem hold for each model. The models among them that exhibit exponential decay of correlations are as follows:
\begin{itemize}
    \item Boolean and Boolean ellipse model: because the points in a Poisson point process are independent;
    \item Mat\'ern cluster model: as there is only local dependency, that is, offspring of different parent points are independent;
    \item Cell model: there is only local dependency within each cell, while points in different cells are independent;
    \item Quermass-interaction processes: if they are in the subcritical regime; that is, if no infinite clusters occur \cite{benes};
    \item DPP model: determinantal point processes with exponentially decreasing kernels exhibit exponential decay of correlations \cite{yukich}, but the Bessel-type kernel used here does not decrease exponentially.
\end{itemize}
Regarding the $(\star)$ condition, it is shown in \cite{tgfpp} that for Mat\'ern cluster and determinantal point processes, the $(\star)$ condition holds. Since the marking is done independently, it also holds for Mat\'ern cluster and $DPP$ models used here. The Boolean, Boolean ellipse, and Cell models satisfy it trivially. To apply Theorem \ref{CLTPD}, it remains to check that the variance condition is satisfied. Our simulation study suggests that, in some cases, even when necessary conditions, such as exponential decay of correlations, are not satisfied, the central limit theorem still holds.

\subsection{Goodness-of-Fit Test Based on Rectangular Partition of PD }\label{partition}
It is of practical use to determine whether there is clustering or repulsion in a random set model. Therefore, our goal in the simulation study is to assess how well the test statistics derived from Theorem \ref{CLTPD} detect whether a realisation of a random set originates from a Boolean model. We also use other models as null models. \\
Assume that we have a realisation of a germ-grain random set model that satisfies the conditions of Theorem \ref{CLTPD}, denoted by $\tilde{\mathcal{P}}_n$, a fixed deterministic radius $r$, and the persistence diagram $PD$ of $\tilde{\mathcal{P}}_n$. We divide the rectangle $[-r,r]^2$ into $m$ rectangles $A_j,j\in \{1,\dots,m\}.$ We then define $m$ functions $f_j,j\in\{1,\dots,m\}$ as $$f_j \colon [-r,r]^2 \to \R, f_j(x,y)=(y-x)\mathbbm{1}_{A_j}(x,y).$$ Following the conclusion of Theorem \ref{CLTPD}, the test statistic $$\sum\limits_{i\in \mathcal I_{q}(\tilde{\mathcal{P}_n})}f_j(b_i,d_i)=\sum\limits_{i\in \mathcal I_{q}(\tilde{\mathcal{P}_n})}(d_i-b_i)\mathbbm{1}_{A_j}(b_i,d_i)$$ is asymptotically normal and therefore
characterised by its mean and variance. Upon standardisation, these random variables converge in distribution to a standard normal variable. As we do not have the theoretical mean and variance of the test statistic for the null model, we estimate them from $1000$ realisations of the Boolean model and from $100$ realisations for other models, and use those estimates $(\mu,\Sigma)$ to test goodness of fit. We test whether the observed data are consistent with the data estimated from the null model. Each observation is centered and standardised using the estimated covariance, and the sum of squared transformed values is compared to a $\chi^2$ distribution with degrees of freedom equal to $m$. More formally, let $X_1, \dots, X_n \in \mathbb{R}^m$ be observed summary statistics, and let $\mu$ and $\Sigma$ be estimates of their mean and covariance under the null hypothesis. We then standardise the observations to obtain $Z_i = (X_i - \mu)\, \Sigma^{-1/2}$. We use the test statistics $T_i = (X_i - \mu)^\top \Sigma^{-1} (X_i - \mu)$, for which $T_i \sim \chi^2_m$ under the null hypothesis.

\begin{figure}[h!]
    \centering
    \includegraphics[scale=0.4]{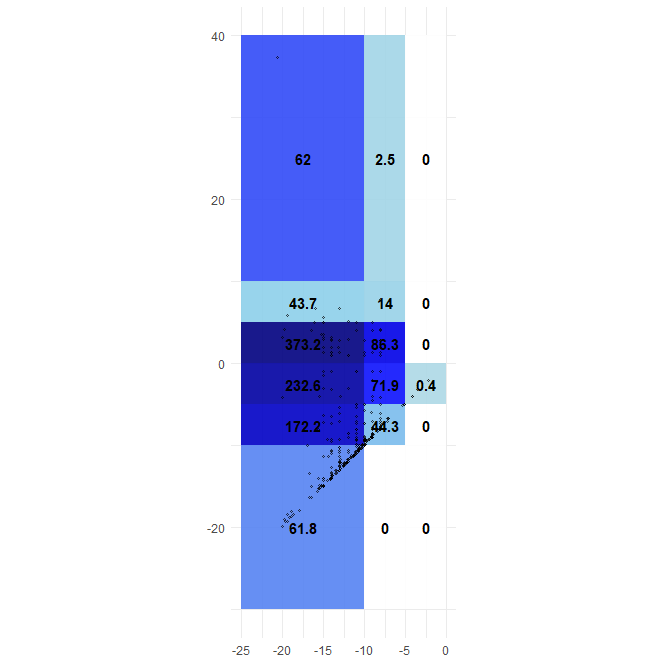}
    \caption{Graphical representation of the division into rectangles of the PD for dimension $0$. The number in the middle of each rectangle shows the estimated expected value of the test statistic for that rectangle, estimated from $1000$ realisations of the Boolean model (the boundary rectangles have been cropped for clearer presentation).}
    \label{PD0_rectangle}
\end{figure}

\begin{figure}[h!]
    \centering
    \includegraphics[scale=0.3]{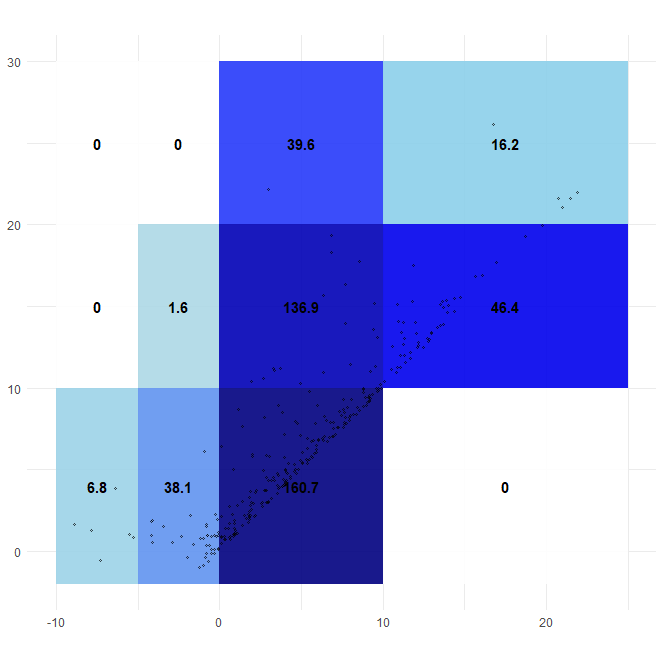}
    \caption{Graphical representation of the division into rectangles of the PD for dimension $1$. The number in the middle of each rectangle shows the estimated expected value of the test statistic for that rectangle, estimated from $1000$ realisations of the Boolean model (the boundary rectangles were cropped for clearer presentation).}
    \label{PD2_rectangle}
\end{figure}

We divide the plane covered with points of $PD$ into rectangles by partitioning the $x$-axis and the $y$-axis in the following way:
\begin{itemize}
    \item Boolean model: $x$-axes at $-40,-10,-5,0$ and $y$-axes at $-30,-10,-5,0,5,10,75$ for dimension $0,$ and $x$-axes at $-10,-5,0,10,50$ and  $y$-axes at $-2,10,20,50$ for dimension $1.$
    \item DPP:  $x$-axes at $-33,-10,-5,0$ and $y$-axes at $-25,-5,0,5,10,50$ for dimension $0,$ and $x$-axes at $-10,-5,0,10,25$ and $y$-axes at $0,10,25$ for dimension $1.$ 
    \item Repulsive: $x$-axes at $-26,-10,-5,0$ and $y$-axes at $-20,-5,0,5,10,30$ for dimension $0,$ and $x$-axes at $-5,0,10,30$ and $y$-axes at $0,10,30$ for dimension $1.$ 
    \item Cluster: $x$-axes at $-30,-10,-5,0$ and $y$-axes at $-30,-5,0,10,50,130$ for dimension $0,$ and $x$-axes at $-10,0,20,50,115$ and $y$-axes at $0,30,70,115$ for dimension $1.$ 
    \item Boolean ellipse: $x$-axes at $-20,-10,-5,0$ and $y$-axes at $-20,-5,0,10,50$ for dimension $0,$ and $x$-axes at $-10,0,10,20,40$ and $y$-axes at $0,10,20,40$ for dimension $1.$ 
    \item Matern Cluster: $x$-axes at $-30,-10,-5,0$ and $y$-axes at $-30,-5,0,20,60$ for dimension $0,$ and $x$-axes at $-15,0,10,20,50$ and $y$-axes at $0,10,20,50$ for dimension $1.$
    \item Cell: $x$-axes at $-25,-10,-5,0$ and $y$-axes at $-25,-5,0,20,55$ for dimension $0,$ and $x$-axes at $-15,0,10,20$ and $y$-axes at $0,5,10,25$ for dimension $1.$
\end{itemize}
Note that the division into rectangles can be done arbitrarily for each model, leading to different results. We divided them as described above after analysing which rectangles might be important for each model when examining the graphical representation of PDs. We omitted the rectangles where the realisations do not contain any points. Figures \ref{PD0_rectangle} and \ref{PD2_rectangle} illustrate how we divided the area covered by the points of PD for dimensions $0$ and $1$, respectively, into rectangles for the Boolean model. The number in the middle of the rectangle shows the estimated expected value of the test statistic for that rectangle, estimated from $1000$ realisations of the Boolean model (we cropped the boundary rectangles for clearer presentation). 

The results of the test using $100$ realisations of each model are shown in Table \ref{rezRecB}. The first column presents the results for the $PD$ for dimension $0$, and the second for dimension $1$. The first number is the percentage of rejection for $p \leq 0.05$, and the second for $p \leq 0.1$, with the Boolean model as the null hypothesis. The results show that this test statistic performs very well in determining whether a realisation comes from the Boolean model. However, it is not perfect for realisations of the Boolean ellipse model for dimension $1$. This is not surprising, as it shares the same underlying point process as the Boolean model, so the test statistic, which is based on detecting the nature of the empty spaces between the components, fails. It also does not recognise the Matern cluster for dimension $1$ very well. Tables \ref{rezRecDPP}, \ref{rezRecBe}, \ref{rezRecMC} and \ref{rezRecCell} show the results of goodness-of-fit tests using different random set models as the null hypothesis in the following order: DPP, Boolean ellipse, Matern cluster, and Cell. When we used the Repulsive and Cluster models as null hypotheses, we obtained a $100\%$ rejection rate for each other random set model. Therefore, we did not include those tables to avoid redundancy. It can be seen that the results are very good, except when comparing Cell to DPP, as they are very similar.

\begin{table}[h!]
\centering
\caption{Table showing results of testing using rectangles (the first number is the percentage of rejection for $p\leq 0.05$ and the second one is for $p \leq 0.1$). (Null hypothesis: Boolean)}
\label{rezRecB}
\label{tab:example}
\begin{tabular}{|l|l|l|}
\hline
\textbf{} & \textbf{dimension $0$} & \textbf{dimension $1$} \\ \hline
DPP & $100\%$ & $100\%$ \\ \hline
Repulsive& $100\%$& $100\%$ \\ \hline
Cluster &  $100\%$& $100\%$\\ \hline
Boolean ellipse & $100\%$ & $62\%,\textcolor{lightgray}{76\%}$ \\ \hline
Matern cluster & $70\%,\textcolor{lightgray}{83\%}$ & $40\%,\textcolor{lightgray}{48\%}$ \\ \hline
Cell & $100\%$ & $100\%$ \\ \hline

\end{tabular}
\end{table}

\begin{table}[h!]
\centering
\caption{Table showing results of testing using rectangles (the first number is the percentage of rejection for $p\leq 0.05$ and the second one is for $p \leq 0.1$). (Null hypothesis: DPP)}
\label{rezRecDPP}
\label{tab:example}
\begin{tabular}{|l|l|l|}
\hline
\textbf{} & \textbf{dimension $0$} & \textbf{dimension $1$} \\ \hline
Boolean & $97\%,\textcolor{lightgray}{98\%}$ & $99\%,\textcolor{lightgray}{99\%}$ \\ \hline
Repulsive& $100\%$& $100\%$ \\ \hline
Cluster &  $100\%$& $100\%$\\ \hline
Boolean ellipse & $100\%$ & $100\%$ \\ \hline
Matern cluster & $100\%$ & $100\%$ \\ \hline
Cell & $32\%,\textcolor{lightgray}{38\%}$ & $35\%,\textcolor{lightgray}{46\%}$ \\ \hline

\end{tabular}
\end{table}

\begin{table}[h!]
\centering
\caption{Table showing results of testing using rectangles (the first number is the percentage of rejection for $p\leq 0.05$ and the second one is for $p \leq 0.1$). (Null hypothesis: Boolean ellipse)}
\label{rezRecBe}
\label{tab:example}
\begin{tabular}{|l|l|l|}
\hline
\textbf{} & \textbf{dimension $0$} & \textbf{dimension $1$} \\ \hline
DPP & $100\%$ & $100\%$  \\ \hline
Boolean & $100\%$ & $65\%,\textcolor{lightgray}{75\%}$ \\ \hline
Repulsive& $100\%$& $100\%$ \\ \hline
Cluster &  $100\%$& $100\%$\\ \hline
Matern cluster & $100\%$ & $97\%,\textcolor{lightgray}{99\%}$  \\ \hline
Cell & $100\%$ & $100\%$ \\ \hline

\end{tabular}
\end{table}

\begin{table}[h!]
\centering
\caption{Table showing results of testing using rectangles (the first number is the percentage of rejection for $p\leq 0.05$ and the second one is for $p \leq 0.1$). (Null hypothesis: Matern Cluster)}
\label{rezRecMC}
\label{tab:example}
\begin{tabular}{|l|l|l|}
\hline
\textbf{} & \textbf{dimension $0$} & \textbf{dimension $1$} \\ \hline
DPP & $100\%$ & $100\%$  \\ \hline
Boolean & $59\%,\textcolor{lightgray}{71\%}$ & $71\%,\textcolor{lightgray}{77\%}$ \\ \hline
Repulsive& $100\%$& $100\%$ \\ \hline
Cluster &  $100\%$& $99\%,\textcolor{lightgray}{99\%}$\\ \hline
Boolean ellipse & $100\%$ & $94\%,\textcolor{lightgray}{97\%}$ \\ \hline
Cell & $100\%$ & $100\%$ \\ \hline

\end{tabular}
\end{table}

\begin{table}[h!]
\centering
\caption{Table showing results of testing using rectangles (the first number is the percentage of rejection for $p\leq 0.05$ and the second one is for $p \leq 0.1$). (Null hypothesis: Cell)}
\label{rezRecCell}
\label{tab:example}
\begin{tabular}{|l|l|l|}
\hline
\textbf{} & \textbf{dimension $0$} & \textbf{dimension $1$} \\ \hline
DPP & $34\%,\textcolor{lightgray}{44\%}$ & $71\%,\textcolor{lightgray}{78\%}$  \\ \hline
Boolean & $90\%,\textcolor{lightgray}{72\%}$ & $100\%$ \\ \hline
Repulsive& $100\%$& $100\%$ \\ \hline
Cluster &  $100\%$& $100\%$\\ \hline
Boolean ellipse & $100\%$ & $100\%$ \\ \hline
Matern cluster & $99\%,\textcolor{lightgray}{100\%}$ & $100\%$  \\ \hline

\end{tabular}
\end{table}

\subsection{Goodness-of-Fit Test Based on Summary Functions of PD}\label{summary}

In our attempt to test the goodness of fit for random set models, we also used different summary functions of $PD$s. Firstly, we consider the accumulated persistence function introduced in \cite{APF}.

For a given persistence diagram $PD^q,$ where $q$ is the dimension of the topological features it captures, we observe each point in the diagram as a triplet $(b_i,d_i,c_i)$, where $c_i$ denotes the multiplicity of each point, $b_i$  denotes  the birth time, and $d_i$  denotes  the death time. We denote by $l_{i} = d_{i} - b_{i}$ the lifetime and by $m_{i} = \frac{b_{i} + d_{i}}{2}$ the mean age of each feature in the diagram, which we refer to as meanage to remain consistent with the terminology in the paper \cite{APF}. Suppose we have $n$ distinct points $(b_i,d_i)$ in the persistence diagram. The accumulated persistence function ($APF$) is defined as \begin{equation}
 APF_q(m)=\sum_{{i}=1}^n c_{i}l_{i}\mathds{1}(m_{i} \leq m), \quad m \in \mathbb R,
\end{equation} where $\mathds{1} ( \cdot)$ is the indicator function, and $q$ denotes the dimension of the topological features under consideration.
$APF_q$ cumulatively sums the lifetimes of the features with respect to their meanage. To use the result of Theorem \ref{CLTPD}, $APF_q$ can be written as
$$APF_q(m)=\sum_{i \in \mathcal{I}_q} f_i(b_i, d_i), \text{ where } f_i = (d_i - b_i)\mathds{1}_{\bigl\{\frac{b_i + d_i}{2} \leq m\bigl\}}(b_i, d_i),$$
where $\mathcal{I}_q$ is an index set over all $q$-dimensional features.

The other summary function is a support function of the lift zonoid of a PD which we introduced in \cite{vesnamarcela}. To obtain it, we first consider the weighted persistence diagram, i.e.
\[
PD^{q}_W(X)=\sum\limits_{i \in \mathcal I_q} (d_i - b_i) \delta_{(b_i, d_i)},\]
to assign greater significance to features that persist longer, as we consider them more important for the shape of a set. The lift zonoid of $P D_W^q(X)$ is defined as: 
 \begin{equation} \label{lz}
 Z=\bigoplus\limits_{i \in \mathcal I_q}(d_i-b_i) \cdot[\mathbf{0},(1,b_i,d_i)].
 \end{equation}
It is the Minkowski sum of the segments in $\mathbb R^3$ with one endpoint at the origin $\mathbf{0}$ and the other at the point $(1,b_i,d_i)$, multiplied by the scalar $l_i=d_i-b_i.$
Its support function $h_Z:S^2\to \mathbb R$ is calculated as follows:
\begin{equation} \label{hlz}
    h_Z^q(u)=\sum\limits_{i \in \mathcal I_q} l_i\max{\{0,\left\langle u,(1,b_i,d_i)\right\rangle\}}, \ u \in S^2,
\end{equation} where $S^2$ stands for the unit sphere. Note that this function is again in a form suitable for using the result of Theorem \ref{CLTPD}.
Since a convex compact set is uniquely determined by its support function, the function $h_Z^q(u)$ defined in (\ref{hlz}) uniquely determines the corresponding lift zonotope in (\ref{lz}). For a detailed analysis of these two summary functions, see \cite{vesnamarcela}. In our simulation study, we parametrise the domain $S^2$ of the support function in the usual way, that is, we take $(\rho, \phi) \in [0,2\pi]\times[0,\pi]$ and identify $u \in S^2$ with $(\rho, \phi)$ such that $u=(\sin(\rho)\cos(\phi),\sin(\rho)\sin(\phi),\cos(\rho)).$

We perform goodness-of-fit testing as in the previous section. We used only the function values for those arguments whose values exhibited approximately normal behaviour, in order to perform the $\chi^2$ test. For $h_Z^0$ we used $\rho=\frac{\pi}{2}, \phi=0,\frac{3\pi}{2}$, for $h_Z^1$ $\rho=\frac{\pi}{2}, \phi=0,\frac{\pi}{4}.$ For $APF_0$ we used $m=-10,30,$ and for $APF_1$ $m=5,40.$ We select those arguments because our graphical analysis indicated that, for these arguments, the function values differ the most.
 Tables \ref{rezB}, \ref{rezDPP}, \ref{rezBe}, \ref{rezMatClust} and \ref{rezCell} show the results of our test using the summary functions mentioned above and $100$ realisations of each random set model, under different null hypotheses. The first number is the percentage of rejection for $p \leq 0.05$, and the second is for $p \leq 0.1$. We omitted tables showing test results when using the repulsive and cluster models as the null model to avoid redundancy, since the percentage of rejection is $100\%$ in all cases except when comparing the Matern cluster to the Cluster model using $APF_1$, where we obtain $88\%,\textcolor{lightgray}{91\%}$. For each random set model, we highlight in bold the highest percentage. For example, when the Boolean model is the null model, $h_Z^1$ outperforms the other summary functions in all cases except for the Boolean ellipse model. This is expected, as in that model the germs are arranged in the same way as in the Boolean model, and $h_Z^1$ detects the nature of the voids between components, so it cannot distinguish the Boolean ellipse model from the Boolean model effectively. 

\begin{table}[h!]
\centering
\caption{Null hypothesis: Boolean}
\label{rezB}
\label{tab:example}
\begin{tabular}{|l|l|l|l|l|}
\hline
\textbf{} & \textbf{$h_Z^0$} & \textbf{$h_Z^1$} & \textbf{$APF_0$} & \textbf{$APF_1$} \\ \hline
DPP & $85\%,\textcolor{lightgray}{94\%}$ & $\mathbf{99\%},\textcolor{lightgray}{100\%}$ & $83\%,\textcolor{lightgray}{87\%}$& $53\%,\textcolor{lightgray}{66\%}$\\ \hline
Boolean & $5\%,\textcolor{lightgray}{9\%}$ & $6\%,\textcolor{lightgray}{10\%}$ & $5\%,\textcolor{lightgray}{8\%}$&  $5\%,\textcolor{lightgray}{9\%}$\\ \hline
Repulsive& $\mathbf{100\%}$& $\mathbf{100\%}$ & $\mathbf{100\%}$ & $\mathbf{100\%}$ \\ \hline
Cluster &  $\mathbf{100\%}$& $\mathbf{100\%}$ & $\mathbf{100\%}$ & $\mathbf{100\%}$\\ \hline
Boolean ellipse & $\mathbf{100\%}$ & $65\%,\textcolor{lightgray}{80\%}$ & $\mathbf{100\%}$ & $55\%,\textcolor{lightgray}{67\%}$\\ \hline
Matern cluster & $47\%,\textcolor{lightgray}{57\%}$ & $\mathbf{89\%},\textcolor{lightgray}{94\%}$ & $84\%,\textcolor{lightgray}{88\%}$& $85\%,\textcolor{lightgray}{91\%}$\\ \hline
Cell & $19\%,\textcolor{lightgray}{39\%}$ & $\mathbf{100\%}$ & $22\%,\textcolor{lightgray}{35\%}$& $43\%,\textcolor{lightgray}{52\%}$\\ \hline

\end{tabular}
\end{table}

\begin{table}[h!]
\centering
\caption{Null hypothesis: DPP}
\label{rezDPP}
\label{tab:example}
\begin{tabular}{|l|l|l|l|l|}
\hline
\textbf{} & \textbf{$h_Z^0$} & \textbf{$h_Z^1$} & \textbf{$APF_0$} & \textbf{$APF_1$} \\ \hline
DPP & $3\%,\textcolor{lightgray}{8\%}$& $5\%,\textcolor{lightgray}{10\%}$& $3\%,\textcolor{lightgray}{9\%}$& $4\%,\textcolor{lightgray}{9\%}$\\ \hline
Boolean & $\mathbf{96\%},\textcolor{lightgray}{97\%}$& $\mathbf{96\%},\textcolor{lightgray}{98\%}$& $80\%,\textcolor{lightgray}{85\%}$&  $78\%,\textcolor{lightgray}{84\%}$\\ \hline
Repulsive& $\mathbf{100\%}$& $\mathbf{100\%}$ & $\mathbf{100\%}$& $\mathbf{100\%}$\\ \hline
Cluster &  $\mathbf{100\%}$& $\mathbf{100\%}$ & $\mathbf{100\%}$& $\mathbf{100\%}$\\ \hline
Boolean ellipse & $\mathbf{100\%}$ & $\mathbf{100\%}$& $\mathbf{100\%}$& $\mathbf{100\%}$\\ \hline
Matern cluster & $\mathbf{100\%}$& $\mathbf{100\%}$& $99\%,\textcolor{lightgray}{100\%}$& $\mathbf{100\%}$\\ \hline
Cell & $\mathbf{53\%},\textcolor{lightgray}{64\%}$& $44\%,\textcolor{lightgray}{60\%}$& $39\%,\textcolor{lightgray}{51\%}$& $27\%,\textcolor{lightgray}{39\%}$\\ \hline

\end{tabular}
\end{table}
\begin{table}[h!]
\centering
\caption{Null hypothesis: Boolean ellipse}
\label{rezBe}
\label{tab:example}
\begin{tabular}{|l|l|l|l|l|}
\hline
\textbf{} & \textbf{$h_Z^0$} & \textbf{$h_Z^1$} & \textbf{$APF_0$} & \textbf{$APF_1$} \\ \hline
DPP & $\mathbf{100\%}$& $\mathbf{100\%}$& $\mathbf{100\%}$& $\mathbf{100\%}$\\ \hline
Boolean & $\mathbf{100\%}$& $78\%,\textcolor{lightgray}{84\%}$& $\mathbf{100\%}$&  $82\%,\textcolor{lightgray}{83\%}$\\ \hline
Repulsive& $\mathbf{100\%}$& $\mathbf{100\%}$ & $\mathbf{100\%}$& $\mathbf{100\%}$\\ \hline
Cluster &  $99\%,\textcolor{lightgray}{99\%}$& $\mathbf{100\%}$ & $\mathbf{100\%}$& $\mathbf{100\%}$\\ \hline
Boolean ellipse & $5\%,\textcolor{lightgray}{11\%}$& $7\%,\textcolor{lightgray}{9\%}$& $5\%,\textcolor{lightgray}{11\%}$& $7\%,\textcolor{lightgray}{11\%}$\\ \hline
Matern cluster & $\mathbf{100\%}$& $99\%,\textcolor{lightgray}{100\%}$& $\mathbf{100\%}$& $96\%,\textcolor{lightgray}{98\%}$\\ \hline
Cell & $\mathbf{100\%}$& $\mathbf{100\%}$& $\mathbf{100\%}$& $\mathbf{100\%}$\\ \hline

\end{tabular}
\end{table}

\begin{table}[h!]
\centering
\caption{Null hypothesis: Matern Cluster}
\label{rezMatClust}
\label{tab:example}
\begin{tabular}{|l|l|l|l|l|}
\hline
\textbf{} & \textbf{$h_Z^0$} & \textbf{$h_Z^1$} & \textbf{$APF_0$} & \textbf{$APF_1$} \\ \hline
DPP & $\mathbf{100\%}$& $\mathbf{100\%}$& $\mathbf{100\%}$& $\mathbf{100\%}$\\ \hline
Boolean & $23\%,\textcolor{lightgray}{35\%}$& $\mathbf{84\%},\textcolor{lightgray}{91\%}$& $59\%,\textcolor{lightgray}{75\%}$&  $82\%,\textcolor{lightgray}{90\%}$\\ \hline
Repulsive& $\mathbf{100\%}$& $\mathbf{100\%}$ & $\mathbf{100\%}$& $\mathbf{100\%}$\\ \hline
Cluster &  $\mathbf{100\%}$& $90\%,\textcolor{lightgray}{95\%}$& $\mathbf{100\%}$& $95\%,\textcolor{lightgray}{98\%}$\\ \hline
Boolean ellipse & $\mathbf{100\%}$& $98\%,\textcolor{lightgray}{100\%}$& $\mathbf{100\%}$& $91\%,\textcolor{lightgray}{94\%}$\\ \hline
Matern cluster & $6\%,\textcolor{lightgray}{12\%}$& $5\%,\textcolor{lightgray}{9\%}$& $5\%,\textcolor{lightgray}{11\%}$& $4\%,\textcolor{lightgray}{10\%}$\\ \hline
Cell & $77\%,\textcolor{lightgray}{90\%}$& $\mathbf{100\%}$& $85\%,\textcolor{lightgray}{93\%}$& $99\%,\textcolor{lightgray}{99\%}$\\ \hline

\end{tabular}
\end{table}

\begin{table}[h!]
\centering
\caption{Null hypothesis: Cell}
\label{rezCell}
\label{tab:example}
\begin{tabular}{|l|l|l|l|l|}
\hline
\textbf{} & \textbf{$h_Z^0$} & \textbf{$h_Z^1$} & \textbf{$APF_0$} & \textbf{$APF_1$} \\ \hline
DPP & $46\%,\textcolor{lightgray}{61\%}$& $\mathbf{72\%},\textcolor{lightgray}{78\%}$& $38\%,\textcolor{lightgray}{48\%}$& $34\%,\textcolor{lightgray}{50\%}$\\ \hline
Boolean & $59\%,\textcolor{lightgray}{73\%}$& $\mathbf{100\%}$& $24\%,\textcolor{lightgray}{39\%}$&  $66\%,\textcolor{lightgray}{74\%}$\\ \hline
Repulsive& $\mathbf{100\%}$& $\mathbf{100\%}$ & $\mathbf{100\%}$& $\mathbf{100\%}$\\ \hline
Cluster &  $\mathbf{100\%}$& $\mathbf{100\%}$& $\mathbf{100\%}$& $\mathbf{100\%}$\\ \hline
Boolean ellipse & $\mathbf{100\%}$& $\mathbf{100\%}$& $\mathbf{100\%}$& $\mathbf{100\%}$\\ \hline
Matern cluster & $96\%,\textcolor{lightgray}{97\%}$& $\mathbf{100\%}$& $91\%,\textcolor{lightgray}{96\%}$& $\mathbf{100\%}$\\ \hline
Cell & $6\%,\textcolor{lightgray}{7\%}$& $5\%,\textcolor{lightgray}{8\%}$& $6\%,\textcolor{lightgray}{12\%}$& $4\%,\textcolor{lightgray}{8\%}$\\ \hline

\end{tabular}
\end{table}

\section{Application to real data concerning histological images of breast tissue}\label{application}

\begin{figure}[H]
    \centering
    \includegraphics[width=\linewidth]{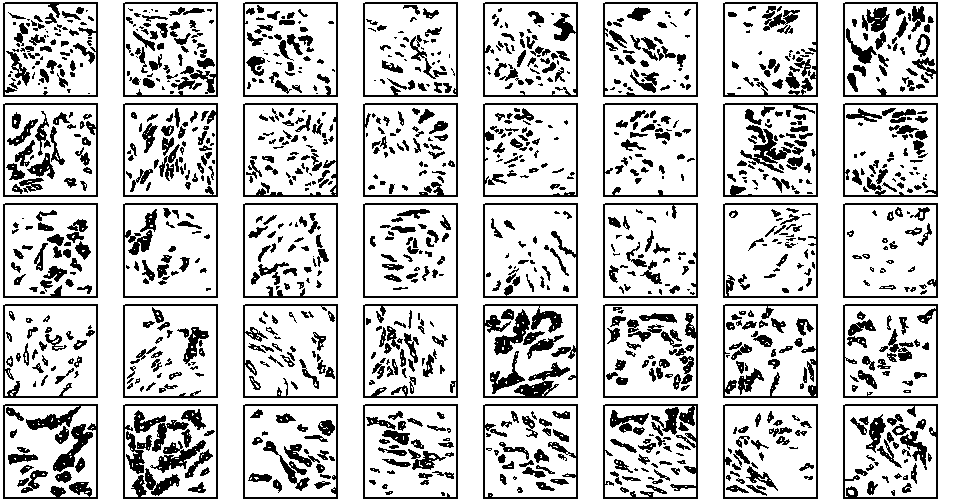}
    \caption{Histological images of mammary breast cancer.}
    \label{fig:mam_realisations}
\end{figure}

\begin{figure}[H]
    \centering
    
        \includegraphics[width=\linewidth]{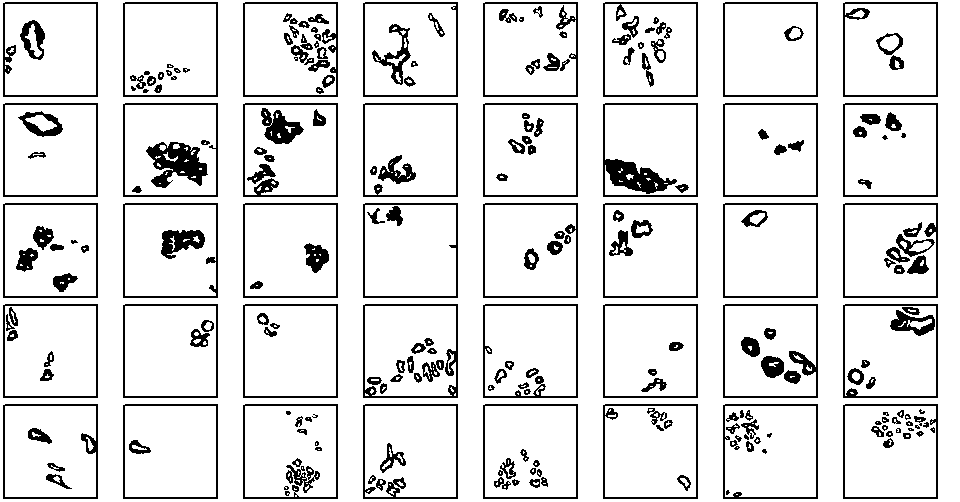}
    \caption{Histological images of mastopathy tissue.}
    \label{fig:mas_realisations}
\end{figure}

We apply our results and methods to real-world data concerning 40 histological images of mammary breast cancer (referred to as Mamca) and 40 images of mastopathy tissue (referred to as Masto), provided by the authors of \cite{mam_mas}.
Since mastopathy is a benign condition, whereas mammary breast cancer represents a malignant transformation of the tissue, distinguishing between these two groups is of substantial clinical importance. Images were converted into $512 \times 512$ binary (black-and-white) pixel matrices, as shown in Figures \ref{fig:mam_realisations} and \ref{fig:mas_realisations}.
 We tested goodness of fit based on a rectangular partition of PD as in Subsection \ref{partition}, but did not obtain good results when partitioning into more than one rectangle. This was because the values do not exhibit approximately normal behaviour, as many rectangles contain no points. Therefore, we examined everything within a single rectangle containing all the points, and in this way, normality was satisfied.
 When testing goodness of fit with Masto as the null model, we obtained rejection rates for Mamca samples $92.5\%,100\%$ for dimension $0,$ and $77.5\%,85\%$ for dimension $1$, for $p\leq 0.05$ and $p\leq 0.1$ respectively. When we conducted the testing with Mamca as the null model, we obtained rejection rates for Masto samples $85\%,90\%$ for dimension $0,$ and $60\%,72\%$ for dimension $1$, for $p\leq 0.05$ and $p\leq 0.1$ respectively. We also attempted to conduct tests using other summary functions as described in Subsection \ref{summary}, but the values of these functions did not exhibit normal behaviour, particularly for Mamca samples. This may be because the observation window is not large enough or due to the nature of these samples, whose distribution is unknown.

\section{Discussion}
In this paper, we established the central limit theorem for functionals of M-bounded persistence diagrams in germ-grain random set models, particularly for those exhibiting exponential decay of correlations. Our work relies on recent results from \cite{prez} on stabilisation methods for marked point processes with exponential decay of correlations.

A simulation study has shown that goodness-of-fit tests relying on asymptotic normality perform well in distinguishing germ-grain models with different interactions. Comparing these results with those of the goodness-of-fit tests in \cite{vesnamarcela}, which use the global envelope tests from \cite{GET}, we observe an overall improvement. Additionally, these new tests are less time-consuming. Our results show that TDA-based statistics offer valuable additional insight for germ-grain random set models exhibiting different interactions.

Finally, we applied tests based on topological data analysis to histological images of breast tissue. These tests showed a good ability to distinguish between mammary breast cancer and mastopathy tissue.

\end{document}